\documentclass[a4paper,review]{amsart}
\usepackage[margin=0.8in]{geometry} 
\newcommand{\eps} {\epsilon}

\usepackage{eurosym}
\usepackage{graphicx}
\usepackage{float}
\usepackage[labelfont=bf]{caption}  

\usepackage[sort]{cite}
\usepackage{subcaption}
\DeclareMathOperator{\supp}{supp}
\usepackage{capt-of}
\usepackage{booktabs}
\usepackage{varwidth}

\usepackage[export]{adjustbox}

\usepackage{tikz}

\newlength\replength

\usepackage{bbm}
\usepackage[utf8]{inputenc}
\usepackage{amssymb}
\usepackage{enumerate}
\usepackage[shortlabels]{enumitem}
\usepackage{mathtools}
\usepackage{graphicx}
\usepackage{url}
\usepackage[utf8]{inputenc}
\usepackage{color}
\definecolor{darkred}{rgb}{0.6,0,0}
\definecolor{darkgreen}{rgb}{0,0.5,0}
\definecolor{darkmagenta}{rgb}{0.5,0,0.5}

\usepackage[pdftex,colorlinks=true,
  linkcolor=blue,
  citecolor=darkgreen,
  urlcolor=darkmagenta]{hyperref}

\makeatletter
\def\@cite#1#2{\textup{[{#1\if@tempswa , #2\fi}]}}
\makeatother

\DeclarePairedDelimiterX\Set[1]\{\}{%
  
  #1
}

\usepackage{xcolor}

\newcommand{\R}{\mathbb R}

\newcommand{\sgn}{\mathop\mathrm{sgn}}

\renewcommand{\d}[1]{\mathinner{\mathrm{d}{#1}}}
\DeclareMathOperator{\TV}{TV}
\DeclareMathOperator{\lipR}{Lip(\R)}

\newtheorem{theorem}{Theorem}[section]
\newtheorem{proposition}[theorem]{Proposition}
\newtheorem{lemma}[theorem]{Lemma}
\newtheorem{definition}[theorem]{Definition}

\newtheorem{remark}[theorem]{Remark}
\newtheorem{corollary}[theorem]{Corollary}

\numberwithin{equation}{section}

\allowdisplaybreaks

\usepackage{tabularx}
\usepackage{multirow}

\begin{document}

\title[Stochastic Nonlocal Conservation Laws with L\'evy Noise]
{Stochastic Nonlocal Conservation Laws with L\'evy Noise: Stability, BV Estimates, and Vanishing Viscosity}

\author[A. Aggarwal]{Aekta Aggarwal}
\address[Aekta Aggarwal]{\newline
Indian Institute of Management Bangalore,
Dorabji Tata School of Undergraduate Studies,
Bengaluru, Karnataka, 560083, India}
\email{aekta.aggarwal@iimb.ac.in}

\author[D. Mukherjee]{Debopriya Mukherjee}
\address[Debopriya Mukherjee]{\newline
Department of Mathematics,
Indian Institute of Technology Indore,
Indore, Madhya Pradesh, 452020, India}
\email{debopriya@iiti.ac.in}

\author[A. A. Panda]{Akash Ashirbad Panda}
\address[Akash Ashirbad Panda]{\newline
Department of Mathematics,
Indian Institute of Technology Indore,
Indore, Madhya Pradesh, 452020, India}
\email{akashpanda@iiti.ac.in}

\date{}




\begin{abstract}
We study one-dimensional stochastic nonlocal conservation laws with a
convolution-type interaction field, driven by multiplicative Brownian noise
and compensated Poisson random measures. The nonlocal operator modifies the
classical entropy structure, while the jump noise produces nonlinear
correction terms in the entropy inequality. We introduce a notion of stochastic entropy solution and establish an $L^1$-stability estimate using spatial doubling of variables and the It\^o--L\'evy formula. This yields uniqueness in the vanishing-viscosity class and, under additional bounds, a localized Cauchy estimate of order $\sqrt{\varepsilon+\varepsilon'}$. Existence is obtained by vanishing viscosity: we construct the viscous solutions through a Faedo--Galerkin scheme and derive viscosity-uniform $L^1$ and localized $BV$ estimates, which allow us to pass to the limit and obtain a spatially $BV$ stochastic entropy solution.
\end{abstract}

\maketitle

\section{Introduction}

Many macroscopic models of transport and collective dynamics are described by
nonlinear hyperbolic conservation laws. A classical example is the
Lighthill--Whitham--Richards (LWR) model for traffic flow,
\begin{eqnarray}\label{LWR}
u_t+(f(u)\nu(u))_x=0,
\quad (t,x)\in Q_T:=(0,T)\times \R,
\end{eqnarray}
where $u$ is the average traffic density, $f\in \lipR$ and
$\nu\in (C^2\cap W^{2,\infty})(\R)$ is the average traffic speed.
The local nature of such models may lead to discontinuities and abrupt
changes in the characteristic velocities. This has motivated conservation
laws with nonlocal fluxes, in which the transport velocity depends on an
averaged state over a spatial interaction region. A prototypical scalar
model is
\begin{eqnarray}
u_t+\left(f(u)\nu(\mu*u)\right)_x&=&0,\label{134}
\end{eqnarray}
where $\nu,\mu\in (C^2\cap W^{2,\infty})(\R)$ and
\[
\mathcal{U}(t,x)=(\mu*u)(t,x)=\int_{\R}\mu(x-\xi)u(t,\xi)\,\d\xi .
\]
The flux at $(t_0,x_0)$ thus depends not only on $u(t_0,x_0)$ but also on
the surrounding state through $\mathcal U(t_0,x_0)$. Such models arise in
traffic with look-ahead interactions, crowd motion, granular flows, laser
technology, sedimentation, conveyor belts, opinion formation, supply chains,
and multilane traffic. Since the nonlocal dependence changes the structure
of the equation, their analysis requires techniques different from the
local theory, and their well-posedness, from scalar equations to systems,
is now well studied via finite-volume approximations, doubling-of-variables
and fixed-point arguments; see
\cite{CG2019,BFK2022,CK24,CCV2022,FGR2021,AHV2025,CGL2012,ACG2015,
BG2016,BHL2023,FGKP2022,AHV2023_1,AHV2024,BBKT2011,CM2015,AS2012,
GHS+2014,ANT2007,AV2023,FCV2023,KP2021,KLS2018,ACT2015,AG2016,AHV2023,
CG2023,AV2026}
and the references therein. More recently, nonlocal conservation laws with
memory effects, where the flux also depends on past states, have been
considered.

Parallel to, and largely independent of, the nonlocal theory described
above, a substantial literature has investigated the effect of random
perturbations on conservation laws. Such perturbations are natural in
transport and collective-motion models, where environmental variability,
unresolved microscopic fluctuations, and abrupt external events are not
captured by the deterministic equation \eqref{LWR}. This has motivated the
study of stochastic conservation laws of the form
\begin{equation}\label{stoch-local}
\d u(t,x)
+\partial_x\!\big(f(u(t,x))\big)\,\d t
=
\Phi(u(t,x))\,\d W(t)
+
\int_{\mathcal Z}
G(u(t,x),z)\,
\widetilde\eta(\d z,\d t).
\end{equation}
Building on Kru\v{z}kov's deterministic entropy formulation
\cite{Kruzkov}, various notions of stochastic entropy solution and their
associated existence and uniqueness theories were developed for
multiplicative Brownian forcing in
\cite{Kim2003,FengNualart2008,DebusscheVovelle2010,Hofmanova2013};
weak-entropy formulations for Cauchy and Dirichlet problems were further
studied in \cite{BVW2012,BVW2014}. A complementary pathwise theory, in
which randomness acts through the flux rather than as a stochastic source,
was developed in \cite{LPS2013,LPS2014}. Conservation laws driven by jump
or L\'evy noise have also been investigated in
\cite{BKM2015,BKM2015dep,BiswasMajee2014,KMV2017,
BiswasKhanVallet2026}.

The distinction between these stochastic mechanisms is important in the
present setting. In \eqref{stoch-local}, and in the model considered below,
the noise acts directly on the state through multiplicative source terms.
Consequently, the entropy balance contains both the usual It\^o correction
associated with the Brownian component and a nonlinear compensator term
generated by the jumps. These contributions are absent from pathwise
random-flux formulations.

The nonlocal and stochastic theories described above have largely developed
along separate directions. This naturally leads to the question of how
entropy admissibility and well-posedness are affected when a
solution-dependent nonlocal transport field and stochastic forcing are
present simultaneously. We therefore consider
\begin{equation}\label{main-eq}
\d u(t,x)
+
\partial_x\!\big(
f(u(t,x))\,\mathcal U(t,x)
\big)\,\d t
=
\Phi(u(t,x))\,\d W(t)
+
\int_{\mathcal Z}
G(u(t,x),z)\,
\widetilde\eta(\d z,\d t),
\end{equation}
where $W$ is a Brownian motion and $\widetilde\eta$ is a compensated
Poisson random measure on the mark space $\mathcal Z$ with compensator
$\nu(\d z)\,\d t$. The nonlocal field $\mathcal U$ represents the spatially
averaged interaction mechanism introduced above, whereas the two stochastic
terms describe qualitatively different forms of uncertainty superimposed on
this transport dynamics.

More precisely, the Brownian term
\[
\Phi(u)\,\d W(t)
\]
models continuous temporal fluctuations whose amplitude depends on the
current state, while
\[
\int_{\mathcal Z}G(u,z)\,\widetilde\eta(\d z,\d t)
\]
accounts for abrupt random events, with the mark $z$ describing their
character or magnitude and the measure $\nu$ determining their statistical
frequency. Since the Brownian motion in \eqref{main-eq} is common to the
spatial domain, it may be interpreted as an aggregate or system-wide
fluctuation with state-dependent local intensity. For instance, in
look-ahead traffic models it may represent persistent variability in driving
behaviour, weather, or road conditions, whereas the jump component may
represent sudden braking events, accidents, lane closures, or external
control actions. Analogously, in crowd dynamics the Brownian component can
represent continuous fluctuations in collective motion, while the jump term
can account for sudden reactions to unexpected events or obstructions.
These interpretations are illustrative; mathematically, the multiplicative
structure expresses that the strength of the perturbation depends on the
current state.

We assume suitable Lipschitz, integrability, and regularity conditions on
the flux $f$, the interaction kernel, and the stochastic coefficients
$\Phi$ and $G$. In particular,
\[
\Phi(0)=0,
\qquad
G(0,z)=0
\quad\text{for }\nu\text{-a.e. }z,
\]
so that the stochastic forcing vanishes at the zero state. No monotonicity
of the noise intensity with respect to the solution is imposed. The precise
probabilistic framework and structural assumptions are stated in
Section~\ref{sec:model}.

Having introduced \eqref{main-eq}, several natural questions arise, which we
answer in this paper.

\begin{enumerate}[(\textbf{Q\arabic*})]

\item \label{q1}
What is the appropriate notion of entropy solution for \eqref{main-eq}, one
that simultaneously accounts for the nonlocal transport operator and the
multiplicative Brownian and L\'evy noise?

\item \label{q2}
Is such a solution unique, and can uniqueness be obtained through a
stochastic analogue of the doubling-of-variables argument, adapted to the
nonlocal flux?

\item \label{q3}
Does \eqref{main-eq} admit a solution at all, and can it be constructed
through a vanishing-viscosity limit of a well-posed viscous approximation?

\item \label{q4}
Can the vanishing-viscosity limit be justified quantitatively, i.e., is
there a computable rate at which the viscous approximations converge?

\end{enumerate}

The analysis developed in this paper answers these questions, combining
techniques from deterministic nonlocal conservation laws, stochastic
analysis, and multiplicative L\'evy-driven SPDEs.

\begin{itemize}

\item We introduce a notion of stochastic entropy solution for
\eqref{main-eq} (Definition~\ref{solution-defn}), built on the shifted entropy flux
$\Psi_k$ rather than the classical Kru\v{z}kov flux difference. This
answers \ref{q1}.

\item We establish an $L^1$-stability estimate (Theorem~\ref{thm:uniq}), yielding
uniqueness within the class of vanishing-viscosity stochastic entropy
solutions. This answers \ref{q2}.

\item We construct viscous Faedo--Galerkin approximations, pass to the
viscous problem, derive viscosity-uniform $L^1$ and localized $BV$ estimates,
and then pass to the vanishing-viscosity limit to obtain a spatially $BV$
stochastic entropy solution (Theorem~\ref{thm:existence}). This answers
\ref{q3}.

\item The stopped comparison estimate yields the localized Cauchy rate
$O(\sqrt{\varepsilon+\varepsilon'})$ for two viscous approximations
(Corollary~\ref{cor:cauchy}); in particular it gives the corresponding
order-$\sqrt{\varepsilon}$ rate when one viscosity parameter is sent to
zero under the stated uniform bounds. This answers \ref{q4}.

\end{itemize}

It should be noted that the present analysis differs from~\cite{LDGW2016,
LGW2021} and~\cite{BGN2025,BGN2026}. The former studied nonlocality arising
from a fractional diffusion operator rather than a convolution in the
transport velocity, so the relevant estimates involve fractional
dissipation and kinetic or fractional-Sobolev structures rather than the
convolution stability estimates central to \eqref{main-eq}. The latter
studied a convolution-flux nonlocal traffic model close in spirit
to~\eqref{134}, with randomness entering through the velocity function
itself, or its parameters, being stochastic (\cite{BGN2025} for Brownian
perturbations, \cite{BGN2026} for more general Markovian ones); the
resulting equations have random transport coefficients but no It\^o
correction or compensated-Poisson terms.

The remainder of the paper is organized as follows.
In Section~\ref{sec:model}, we introduce the mathematical model, state the
standing assumptions, and present the probabilistic and functional framework
used throughout the paper.
Section~\ref{sec:entropy} introduces the notion of stochastic entropy
solutions and derives the corresponding entropy formulation.
The uniqueness theorem is established in Section~\ref{sec:uniq} by means of a
stochastic doubling-of-variables argument, with the auxiliary estimates it
relies on collected and proved separately in Section~\ref{sec:lemmas}.
Section~\ref{sec:existence} is devoted to the existence theory.
We first study the viscous approximation through a Galerkin scheme, derive
uniform estimates independent of both the Galerkin dimension and the
viscosity parameter, establish compactness of the approximate solutions, and
finally pass to the vanishing-viscosity limit to obtain stochastic entropy
solutions.
Section~\ref{sec:conclusion} closes the paper with concluding remarks and
directions for future work.

\smallskip

\section{Mathematical Formulation and Assumptions}\label{sec:model}

Let $\mathfrak P :=
(\Omega,\mathcal F,\mathbb F=\{\mathcal F_t\}_{t\ge0},\mathbb P)$
be a filtered probability space satisfying the usual conditions.
Let $W=\{W(t)\}_{t\ge0}$ be a one-dimensional $\mathbb F$-Brownian
motion. Let $(Z,\mathcal B(Z))$ be a locally compact Polish space and
let $\eta(\d z,\d t)$ be an $\mathbb F$-Poisson random measure on
$Z\times(0,\infty)$ with compensator $\nu(\d z)\d t$, where $\nu$ is
$\sigma$-finite. We write
\[
\widetilde\eta(\d z,\d t)
:=
\eta(\d z,\d t)-\nu(\d z)\d t
\]
for the associated compensated Poisson random measure.

We consider the stochastic nonlocal conservation law
\begin{align}
\label{eq:u}
\d u(t,x)
+
\partial_x\!\Big(
f(u(t,x))\,\mathcal U[u](t,x)
\Big)\,\d t
&=
\Phi(u(t,x))\,\d W(t)
+
\int_Z
G(u(t,x),z)\,
\widetilde\eta(\d z,\d t),
\\
u(0,x)
&=
u_0(x),
\qquad x\in\R .
\label{eq:id}
\end{align}
Here $f:\R\to\R$ is the flux function,
$\Phi:\R\to\R$ is the Brownian coefficient,
$G:\R\times Z\to\R$ is the jump amplitude, and
$\mathcal U[u]$ denotes the nonlocal velocity generated by the solution.
The bracket notation emphasizes the dependence of $\mathcal U$ on $u$.

Throughout the paper, $C>0$ denotes a generic constant whose value may
change from line to line.

\subsection{Nonlocal velocity field and kernel assumptions}
\label{subsec:kernel}

The nonlocal field is given by
\begin{equation}\label{eq:Uconv}
\mathcal U[u](t,x)
:=
(\mu*u(t,\cdot))(x)
=
\int_{\R}
\mu(x-y)u(t,y)\,\d y,
\qquad
(t,x)\in[0,T]\times\R,
\end{equation}
where $\mu:\R\to\R$ is a deterministic interaction kernel.

We impose
\[
\begin{array}{ll}
\textbf{(K1)}\quad \mu\in W^{1,1}(\R),
&
\textbf{(K2)}\quad \mu\in L^\infty(\R),
\\[1mm]
\textbf{(K3)}\quad \mu'\in L^\infty(\R),
&
\textbf{(K4)}\quad \mu''\in L^1(\R).
\end{array}
\]
Assumptions \textbf{(K1)}--\textbf{(K3)} are used in the basic
well-posedness and stability analysis, whereas \textbf{(K4)} is needed
only in estimates involving $\partial_{xx}\mathcal U$, in particular
the $BV$-estimate of Proposition~\ref{prop:BV}.

Since \textbf{(K1)}--\textbf{(K3)} imply
$\mu,\mu'\in L^1(\R)\cap L^\infty(\R)$, we also have
\begin{equation}\label{eq:K-L2}
\|\mu\|_{L^2}^2
\le
\|\mu\|_{L^1}\|\mu\|_{L^\infty},
\qquad
\|\mu'\|_{L^2}^2
\le
\|\mu'\|_{L^1}\|\mu'\|_{L^\infty}.
\end{equation}
these consequences will be used in the truncated Galerkin construction.

For $u,v\in L^1(\R)$, Young's inequality yields
\begin{align}
\|\mathcal U[u]\|_{L^\infty}
&\le
\|\mu\|_{L^\infty}\|u\|_{L^1},
\label{eq:U-Linfty-L1}
\\
\|\partial_x\mathcal U[u]\|_{L^\infty}
&\le
\|\mu'\|_{L^\infty}\|u\|_{L^1},
\label{eq:Ux-Linfty}
\\
\|\mathcal U[u]-\mathcal U[v]\|_{L^\infty}
&\le
\|\mu\|_{L^\infty}\|u-v\|_{L^1},
\label{eq:U-Linfty-stab}
\\
\|\partial_x\mathcal U[u]
-\partial_x\mathcal U[v]\|_{L^\infty}
&\le
\|\mu'\|_{L^\infty}\|u-v\|_{L^1},
\label{eq:Ux-Linfty-stab}
\\
\|\partial_x\mathcal U[u]
-\partial_x\mathcal U[v]\|_{L^1}
&\le
\|\mu'\|_{L^1}\|u-v\|_{L^1},
\label{eq:Ux-L1-stab-basic}
\\
\|\partial_{xx}\mathcal U[u]\|_{L^1}
&\le
\|\mu''\|_{L^1}\|u\|_{L^1}.
\label{eq:Uxx-L1}
\end{align}
The last estimate requires \textbf{(K4)}. These bounds will be used
repeatedly without further comment.

\subsection{Standing assumptions}\label{subsec:assumptions}

We impose the following hypotheses.

\begin{enumerate}[label=\textbf{(A\arabic*)}, leftmargin=2.4em]

\item\label{A1}
\textbf{Flux.}
The function $f:\R\to\R$ is globally Lipschitz, satisfies $f(0)=0$, and
\[
|f(r)-f(s)|
\le
L_f|r-s|,
\qquad r,s\in\R .
\]
By Rademacher's theorem, $f$ is differentiable a.e. on $\R$ and
\[
|f'(r)|\le L_f
\quad\text{for a.e. }r\in\R,
\qquad
|f(r)|\le L_f|r|.
\]
Here and below, $f'$ denotes this a.e. derivative. No boundedness assumption
on $f$ itself is imposed.

\item\label{A2}
\textbf{Brownian coefficient.}
The function $\Phi\in C^1(\R)$ satisfies $\Phi(0)=0$ and
\[
|\Phi(r)-\Phi(s)|
\le
L_\Phi|r-s|,
\qquad r,s\in\R .
\]
Consequently,
\[
|\Phi(r)|\le L_\Phi|r|,
\qquad
\|\Phi'\|_{L^\infty}\le L_\Phi .
\]

\item\label{A3}
\textbf{Jump coefficient.}
The mapping $G:\R\times Z\to\R$ is measurable,
$G(0,z)=0$ for $\nu$-a.e.\ $z$, and
\begin{align*}
\int_Z
|G(r,z)-G(s,z)|^2\,\nu(\d z)
&\le
L_G^2|r-s|^2,
\\
\int_Z
|G(r,z)-G(s,z)|\,\nu(\d z)
&\le
\widetilde L_G|r-s|,
\qquad r,s\in\R .
\end{align*}
In particular,
\[
\int_Z|G(r,z)|^2\,\nu(\d z)
\le
L_G^2|r|^2,
\qquad
\int_Z|G(r,z)|\,\nu(\d z)
\le
\widetilde L_G|r|.
\]

\item\label{A4}
\textbf{Nonlocal field.}
The velocity field is given by
$\mathcal U[u]=\mu*u$, where $\mu$ satisfies
\textbf{(K1)}--\textbf{(K3)}.
Assumption \textbf{(K4)} is imposed whenever the $BV$-estimate is used.

\item\label{A5}
\textbf{Initial datum.}
The random variable $u_0$ is $\mathcal F_0$-measurable and
\[
u_0
\in
L^2\!\big(
\Omega;
L^1(\R)\cap L^2(\R)\cap BV(\R)
\big).
\]

\end{enumerate}


\smallskip
\noindent
\textbf{Assumption used in the $BV$-analysis.}
Whenever the spatial $BV$-estimate, or a subsequent argument relying on it,
is invoked, we additionally impose the kernel condition \textbf{(K4)}. No
$C^2$-regularity of $f$ and no differentiability of $G$ with respect to the
state variable are required: the $BV$ argument below is carried out by
spatial translations and uses only the Lipschitz assumptions
\textbf{(A1)}--\textbf{(A3)}, together with the $BV$-regularity of $u_0$
already included in \textbf{(A5)}.

For $w\in L^1_{\mathrm{loc}}(\R)$, we write
\[
\TV(w)
:=
\sup\left\{
\int_\R w(x)\varphi'(x)\,\d x:
\varphi\in C_c^1(\R),\
\|\varphi\|_{L^\infty}\le1
\right\}.
\]


\smallskip

\section{Formulation of Stochastic Entropy Solutions}
\label{sec:entropy}

We now introduce the entropy formulation associated with
\eqref{eq:u}--\eqref{eq:id}. The main difference from the classical
scalar conservation law is that the flux
\[
F(t,x,r)
=
f(r)\mathcal U[u](t,x)
\]
depends explicitly on the spatial variable through the nonlocal velocity
field $\mathcal U[u]$. Consequently, the usual entropy flux alone does
not account for the entire transport contribution: an additional
zeroth-order term involving $\partial_x\mathcal U[u]$ naturally appears
in the entropy balance. We first derive this structure for smooth
solutions and then use it to define stochastic entropy solutions.

Let $Q_T:=(0,T)\times\R .$ We denote by
$\mathcal N_\Omega^2(0,T;L^1(\R))$ the space of predictable processes $u:\Omega\times(0,T)\longrightarrow L^1(\R)$ such that
\[
\mathbb E
\int_0^T
\|u(t)\|_{L^1(\R)}^2\,\d t
<\infty .
\]

\subsection{Entropy--flux pairs}

For the moment, regard $\mathcal U=\mathcal U(t,x)$ as a given
space--time dependent coefficient. Its dependence on the solution will
be reinstated when the entropy formulation is applied to
$\mathcal U[u]$.

\begin{definition}[\textbf{Entropy--flux pair}]
\label{def:efp}
Let $\beta\in C^2(\R)$ be convex and consider the flux $F(t,x,r)
:= f(r)\mathcal U(t,x).$ The entropy flux associated with $\beta$ is defined by
\begin{equation}\label{eq:qdef}
q(t,x,r)
:=
\int_0^r
\beta'(s)\partial_sF(t,x,s)\,\d s
=
\mathcal U(t,x)
\int_0^r
\beta'(s)f'(s)\,\d s .
\end{equation}
Thus, with
\[
Q_\beta(r)
:=
\int_0^r
\beta'(s)f'(s)\,\d s,
\]
we have
\[
q(t,x,r)
=
\mathcal U(t,x)Q_\beta(r),
\]
and hence, with $r$ kept fixed in the $x$-derivative,
\[
q_r(t,x,r)
=
\beta'(r)f'(r)\mathcal U(t,x),
\qquad
q_x(t,x,r)
=
\partial_x\mathcal U(t,x)Q_\beta(r).
\]

For the shifted entropy
\[
r\longmapsto\beta(r-k),
\qquad k\in\R,
\]
we introduce the shifted entropy flux
\begin{equation}\label{eq:Psi-def}
\Psi_k(r)
:=
\int_k^r
\beta'(s-k)f'(s)\,\d s .
\end{equation}
Accordingly,
\[
q_k(t,x,r)
:=
\mathcal U(t,x)\Psi_k(r).
\]
The lower limit $k$ is convenient because it gives $\Psi_k(k)=0$ and makes the limiting Kru\v{z}kov structure transparent.

In general,
\[
\Psi_k(r)
\neq
\beta'(r-k)\big(f(r)-f(k)\big),
\]
since $\beta'$ varies along the integration interval. For the
Kru\v{z}kov entropy $\beta(r)=|r|$, however,
\[
\Psi_k(r)
=
\sgn(r-k)\big(f(r)-f(k)\big).
\]
Moreover,
\begin{equation}\label{eq:Psi-bounds}
|\Psi_k(r)|
\le
\|\beta'\|_{L^\infty}L_f|r-k|,
\end{equation}
and, for $\beta$ in the admissible class introduced below,
\begin{equation}\label{eq:Psi-correction-bound}
\big|
\beta'(r-k)f(r)-\Psi_k(r)
\big|
\le
C_\beta
\big(
|f(k)|+L_f|r-k|
\big),
\end{equation}
where $C_\beta$ depends only on $\|\beta'\|_{L^\infty}$.
In particular, all transport terms appearing in the entropy formulation
have at most linear growth in $|r-k|$.

For the smooth approximations $\beta_\delta$ of the Kru\v{z}kov
entropy used in Section~\ref{sec:uniq},
\begin{equation}\label{eq:Psi-Kruzkov-limit}
\Psi_k(r)
\longrightarrow
\sgn(r-k)
\big(f(r)-f(k)\big)
\qquad
\text{as }\delta\downarrow0 .
\end{equation}
\end{definition}

\begin{remark}[Admissible entropies]
\label{rem:admissible-beta}
Throughout the entropy formulation we use the class
\[
\mathcal E
:=
\left\{
\beta\in C^2(\R):
\begin{array}{l}
\beta\ge0,\quad
\beta\ \text{convex},\quad
\beta(0)=0,\\[1mm]
\beta'\in L^\infty(\R),\quad
(1+|r|)\beta''(r)\in L^\infty(\R)
\end{array}
\right\}.
\]
Since $\beta$ is convex, nonnegative, and has its minimum at the
origin, we also have $\beta'(0)=0.$ The boundedness of $\beta'$ implies the linear-growth estimate
\[
0\le\beta(r)
\le
\|\beta'\|_{L^\infty}|r|.
\]

The weighted condition on $\beta''$ is imposed in order to control the
It\^o correction under the natural $L^1$-integrability of the entropy
solution. Indeed, since $\Phi(0)=0$ and $|\Phi(r)|\le L_\Phi|r|,$ for every fixed $k\in\R$,
\begin{equation}\label{eq:ito-entropy-integrability}
\beta''(r-k)\Phi(r)^2
\le
C_{\beta,k}\big(1+|r|\big).
\end{equation}
Hence the Brownian It\^o correction is locally integrable in space for
$L^1$-valued solutions.

Similarly, the convexity of $\beta$ and
$|\beta'|\le\|\beta'\|_{L^\infty}$ imply
\[
0
\le
\beta(a+b)-\beta(a)-\beta'(a)b
\le
2\|\beta'\|_{L^\infty}|b|,
\]
so that the jump compensator is controlled by the
$L^1(Z,\nu)$-growth assumption in \textbf{(A3)}.

The class $\mathcal E$ contains, in particular, the regularized
Kru\v{z}kov entropies
\[
\beta_\delta(r)
=
\sqrt{r^2+\delta^2}-\delta
\]
used in the uniqueness analysis.
\end{remark}

\begin{remark}[Temporal regularity]
\label{rem:path-reg}
Because of the Poisson noise, the solution is not expected to have
continuous trajectories. The compactness construction of
Section~\ref{sec:existence} naturally yields c\`adl\`ag paths in the
Hermite--Sobolev space $\mathcal S_{-1}$ introduced in
\eqref{eq:hermite-scale}. We therefore impose
\[
u\in D([0,T];\mathcal S_{-1})
\qquad\mathbb P\text{-a.s.},
\]
rather than the stronger requirement
$u\in D([0,T];L^1(\R))$.

The space $\mathcal S_{-1}$ is sufficiently weak to accommodate the
compactness argument, while still allowing the weak formulation and
the temporal traces required below.
\end{remark}

\subsection{Entropy balance for the viscous approximation}

To motivate the entropy inequality, consider the viscous approximation
\begin{align}
\label{eq:ueps}
\d u^\eps
+
\partial_x
\Big(
f(u^\eps)\mathcal U[u^\eps]
\Big)\,\d t
&=
\eps\partial_{xx}u^\eps\,\d t
+
\Phi(u^\eps)\,\d W(t)
+
\int_Z
G(u^\eps_-,z)\,
\widetilde\eta(\d z,\d t),
\\
u^\eps(0,x)
&=
u_0(x).
\nonumber
\end{align}
Here and below,
\[
u^\eps_-(t,x):=u^\eps(t-,x)
\]
denotes the left limit, which is the predictable quantity appearing in
the Poisson stochastic integral.

For notational convenience, set
\[
\mathcal U^\eps
:=
\mathcal U[u^\eps] \qquad \mbox{and} \qquad Q_\beta(r)
:=
\int_0^r
\beta'(s)f'(s)\,\d s.
\]
Then
\[
q(t,x,r)
=
\mathcal U^\eps(t,x)Q_\beta(r).
\]

The essential point is the following chain-rule identity. For a smooth
function $u$,
\begin{align}
\beta'(u)
\partial_x
\big(
f(u)\mathcal U
\big)
&=
\beta'(u)f'(u)\mathcal U\,\partial_xu
+
\beta'(u)f(u)\partial_x\mathcal U
\nonumber\\
&=
\partial_x
\big(
\mathcal UQ_\beta(u)
\big)
+
\big[
\beta'(u)f(u)-Q_\beta(u)
\big]
\partial_x\mathcal U .
\label{eq:chain-visc}
\end{align}
The second term on the right-hand side of
\eqref{eq:chain-visc} is absent for a spatially homogeneous flux and is
precisely the correction generated by the spatial dependence of the
nonlocal velocity.

Applying the It\^o--L\'evy formula to $\beta(u^\eps)$ and using
\eqref{eq:chain-visc}, we obtain formally
\begin{align}
\label{eq:ito-beta}
\d\beta(u^\eps)
&+
\partial_x
\big(
\mathcal U^\eps Q_\beta(u^\eps)
\big)\,\d t
+
\big[
\beta'(u^\eps)f(u^\eps)-Q_\beta(u^\eps)
\big]
\partial_x\mathcal U^\eps\,\d t
\nonumber\\
&=
\eps\partial_{xx}\beta(u^\eps)\,\d t
-
\eps\beta''(u^\eps)
|\partial_xu^\eps|^2\,\d t
\nonumber\\
&\quad
+
\beta'(u^\eps)\Phi(u^\eps)\,\d W(t)
+
\frac12
\beta''(u^\eps)\Phi(u^\eps)^2\,\d t
\nonumber\\
&\quad
+
\int_Z
\Big[
\beta\big(
u^\eps_-+G(u^\eps_-,z)
\big)
-
\beta(u^\eps_-)
\Big]
\widetilde\eta(\d z,\d t)
\nonumber\\
&\quad+
\int_Z
\Big[
\beta\big(
u^\eps_-+G(u^\eps_-,z)
\big)
-
\beta(u^\eps_-)
-
\beta'(u^\eps_-)G(u^\eps_-,z)
\Big]
\nu(\d z)\,\d t .
\end{align}

For the entropy centered at an arbitrary level $k\in\R$, it is more
convenient to use the shifted flux \eqref{eq:Psi-def} directly.
Since $\frac{\mathrm{d}}{\mathrm{d}r}\Psi_k(r)
= \beta'(r-k)f'(r)\,,$ the analogue of \eqref{eq:chain-visc} is
\begin{equation}\label{eq:shifted-chain}
\beta'(u-k)
\partial_x
\big(
f(u)\mathcal U
\big)
=
\partial_x
\big(
\mathcal U\Psi_k(u)
\big)
+
\big[
\beta'(u-k)f(u)-\Psi_k(u)
\big]
\partial_x\mathcal U .
\end{equation}
Consequently,
\begin{align}
\label{eq:ito-beta-shift}
\d\beta(u^\eps-k)
&+
\partial_x
\Big(
\mathcal U^\eps\Psi_k(u^\eps)
\Big)\,\d t +
\big[
\beta'(u^\eps-k)f(u^\eps)
-
\Psi_k(u^\eps)
\big]
\partial_x\mathcal U^\eps\,\d t
=
\mathcal R_k^\eps ,
\end{align}
where $\mathcal R_k^\eps$ denotes the viscosity, Brownian, and
It\^o--L\'evy contributions on the right-hand side of
\eqref{eq:ito-beta}, with $\beta$, $\beta'$, and $\beta''$ evaluated at
the shifted argument.

The term
\[
-\eps
\beta''(u^\eps-k)
|\partial_xu^\eps|^2
\]
is nonpositive because $\beta$ is convex. After multiplication by a
nonnegative test function and integration by parts, this entropy
dissipation has the favorable sign and gives rise to an inequality in
the limit $\eps\downarrow0$. This motivates the definition below.

\subsection{Stochastic entropy solutions}

\begin{definition}[\textbf{Stochastic entropy solution}]
\label{solution-defn}
Assume \textbf{(A1)}--\textbf{(A5)}.
An $\mathbb F$-adapted process $u$ satisfying
\[
u
\in
\mathcal N_\Omega^2(0,T;L^1(\R))
\cap
L^2
\big(
\Omega;L^\infty(0,T;L^1(\R))
\big), \qquad \mbox{and} \qquad u
\in
D([0,T];\mathcal S_{-1})
\qquad
\mathbb P\text{-a.s.},
\]
is called a \emph{stochastic entropy solution} of
\eqref{eq:u}--\eqref{eq:id} if, for every
\[
\beta\in\mathcal E,
\qquad
k\in\R,
\qquad
0\le\phi\in C_c^\infty([0,T)\times\R),
\]
the following inequality holds $\mathbb P$-a.s.:
\begin{align}
\label{eq:entropy-ineq}
&\int_\R
\beta(u_0-k)\phi(0,x)\,\d x
+
\int_{Q_T}
\beta(u-k)
\partial_t\phi\,\d x\,\d t
\nonumber\\
&\quad
+
\int_{Q_T}
\mathcal U[u]\Psi_k(u)
\partial_x\phi\,\d x\,\d t
\nonumber\\
&\quad
-
\int_{Q_T}
\Big[
\beta'(u-k)f(u)-\Psi_k(u)
\Big]
\partial_x\mathcal U[u]\,
\phi\,\d x\,\d t
\nonumber\\
&\quad
+
\int_0^T
\left(
\int_\R
\beta'(u-k)
\Phi(u)\phi\,\d x
\right)
\d W(t)
\nonumber\\
&\quad
+
\frac12
\int_{Q_T}
\beta''(u-k)
\Phi(u)^2
\phi\,\d x\,\d t
\nonumber\\
&\quad
+
\int_0^T\!\!\int_Z
\left(
\int_\R
\Big[
\beta\big(
u_-(t,x)-k+G(u_-(t,x),z)
\big)
-
\beta\big(
u_-(t,x)-k
\big)
\Big]
\phi(t,x)\,\d x
\right)
\widetilde\eta(\d z,\d t)
\nonumber\\
&\quad
+
\int_0^T\!\!\int_Z\!\!\int_\R
\Big[
\beta\big(
u_-(t,x)-k+G(u_-(t,x),z)
\big)
-
\beta\big(
u_-(t,x)-k
\big)
\nonumber\\
&\hspace{5.0cm}
-
\beta'\big(
u_-(t,x)-k
\big)
G(u_-(t,x),z)
\Big]
\phi(t,x)
\,\d x\,\nu(\d z)\,\d t
\ge0 .
\end{align}
Here $\Psi_k$ is defined in \eqref{eq:Psi-def}.
\end{definition}

\begin{remark}[Well-definedness of the entropy terms]
\label{rem:entropy-well-defined}
The integrability assumptions in Definition~\ref{solution-defn} are
chosen so that every term in \eqref{eq:entropy-ineq} is well defined.

Indeed, since $\phi$ is compactly supported in space and
$\beta'$ is bounded, the Brownian integrand satisfies
\[
\left|
\int_\R
\beta'(u-k)\Phi(u)\phi\,\d x
\right|
\le
C_{\beta,\phi}
\int_{\supp{\phi}}
|u|\,\d x ,
\]
and therefore belongs to $L^2(\Omega\times(0,T))$ by the definition of
$\mathcal N_\Omega^2(0,T;L^1(\R))$.

Likewise, \eqref{eq:ito-entropy-integrability} controls the Brownian
It\^o correction, whereas
\[
\big|
\beta(a+b)-\beta(a)
\big|
\le
\|\beta'\|_{L^\infty}|b| \qquad \mbox{and} \qquad 0
\le
\beta(a+b)-\beta(a)-\beta'(a)b
\le
2\|\beta'\|_{L^\infty}|b|
\]
show that the two jump contributions are integrable under
\textbf{(A3)}. Finally, the bounds
\eqref{eq:Psi-bounds}--\eqref{eq:Psi-correction-bound}, together with
the convolution estimates of Section~\ref{subsec:kernel}, control the
two nonlocal transport terms.
\end{remark}

\begin{remark}[Structure of the nonlocal entropy flux]
\label{rem:flux-correction}
To clarify the role of the additional transport term, consider first a
smooth solution of the deterministic equation
\[
\partial_tu
+
\partial_x
\big(
f(u)\mathcal U[u]
\big)
=
0.
\]
Using \eqref{eq:shifted-chain}, one obtains the exact entropy identity
\begin{equation}\label{eq:entropy-identity}
\partial_t\beta(u-k)
+
\partial_x
\big(
\mathcal U[u]\Psi_k(u)
\big)
+
\Big[
\beta'(u-k)f(u)-\Psi_k(u)
\Big]
\partial_x\mathcal U[u]
=
0 .
\end{equation}
Thus the term
\[
\Big[
\beta'(u-k)f(u)-\Psi_k(u)
\Big]
\partial_x\mathcal U[u]
\]
is not an auxiliary correction introduced for technical reasons; it is
an intrinsic part of the entropy balance for the spatially varying
nonlocal flux.

For the Kru\v{z}kov entropy $\beta(r)=|r|$,
\[
\Psi_k(u)
=
\sgn(u-k)
\big(
f(u)-f(k)
\big),
\]
and therefore
\[
\beta'(u-k)f(u)-\Psi_k(u)
=
\sgn(u-k)f(k).
\]
Accordingly, the two transport terms in the weak entropy formulation
reduce to
\[
\mathcal U[u]\,
\sgn(u-k)
\big(
f(u)-f(k)
\big)
\partial_x\phi
-
\sgn(u-k)f(k)
\partial_x\mathcal U[u]\,
\phi.
\]
This is precisely the nonlocal analogue of the usual Kru\v{z}kov
entropy structure.
\end{remark}

\begin{remark}[Relation with the viscous approximation]
\label{rem:entropy-viscosity}
For smooth viscous solutions, \eqref{eq:entropy-ineq} follows from
\eqref{eq:ito-beta-shift} after multiplying by a nonnegative test
function and integrating over space and time. The viscosity generates
the additional dissipation
\[
\eps
\int_{Q_T}
\beta''(u^\eps-k)
|\partial_xu^\eps|^2
\phi\,\d x\,\d t
\ge0,
\]
which is discarded when passing to the entropy inequality.
Consequently, the formulation in
Definition~\ref{solution-defn} is consistent with the
vanishing-viscosity construction carried out in
Section~\ref{sec:existence}.
\end{remark}

The entropy formulation above therefore incorporates simultaneously
the solution-dependent nonlocal flux, the Brownian It\^o correction,
and the L\'evy jump contribution. It provides the natural framework
for the $L^1$-stability and uniqueness analysis developed in
Section~\ref{sec:uniq}.

\medskip

\section{Uniqueness and $L^1$-Stability}
\label{sec:uniq}

In this section we establish an $L^1$-stability estimate and, as a
consequence, uniqueness in the class of vanishing-viscosity solutions.
The main difficulty is caused by the fact that the flux
\[
f(u)\mathcal U[u]
\]
depends nonlocally on the solution itself. Thus, when two solutions are
compared, one must control not only the difference of their local
states, but also the difference of the associated nonlocal velocity
fields.

A second difficulty is stochastic. A direct doubling-of-variables
argument in both space and time would introduce random entropy levels
evaluated at different times and would lead to an adaptedness problem.
We avoid this by doubling only the spatial variable and comparing the
two viscous approximations at the same time. More precisely, we apply
the It\^o--L\'evy formula directly to
\[
(u^\varepsilon(t,x),v^{\varepsilon'}(t,y))
\longmapsto
\beta_\delta
\big(
u^\varepsilon(t,x)-v^{\varepsilon'}(t,y)
\big).
\]
This equal-time comparison is particularly important for the
stochastic terms: the Brownian correction naturally contains
\[
\big(
\Phi(u^\varepsilon)-\Phi(v^{\varepsilon'})
\big)^2,
\]
while the jump contribution contains
\[
G(u^\varepsilon,z)-G(v^{\varepsilon'},z),
\]
which are precisely the difference structures needed for an
$L^1$-estimate.

Since this argument uses the equations satisfied by the viscous
approximations rather than merely the limiting entropy inequalities,
uniqueness is first proved in the vanishing-viscosity class.

\subsection{Vanishing-viscosity solutions and the stability theorem}

\begin{definition}[Vanishing-viscosity solution]
\label{def:vv}
Let $u_0$ satisfy \textbf{(A5)}. A process $u$ is called a
\emph{vanishing-viscosity solution} of
\eqref{eq:u}--\eqref{eq:id} if there exist a sequence
$\varepsilon_n\downarrow0$ and corresponding $\mathcal F_0$-measurable
regularizations $u_0^{\varepsilon_n}$ such that
\[
u_0^{\varepsilon_n}
\longrightarrow
u_0
\qquad\text{in }
L^2\!\big(
\Omega;L^1(\R)\cap L^2(\R)
\big),
\]
and
\[
\TV(u_0^{\varepsilon_n})
\le
\TV(u_0)
\qquad\text{a.s. for every }n,
\]
while the solutions $u^{\varepsilon_n}$ of the corresponding viscous
problems satisfy
\[
u^{\varepsilon_n}
\longrightarrow
u
\qquad
\text{in }
L^1\!\big(
\Omega\times(0,T)\times I
\big)
\]
for every bounded interval $I\subset\R$.

A canonical choice is the Friedrichs regularization
$u_0^\varepsilon=\rho_\varepsilon*u_0$, for which the above variation
bound holds. By the construction of Section~\ref{subsec:eps-limit}, every
vanishing-viscosity solution is a stochastic entropy solution in the sense
of Definition~\ref{solution-defn}.
\end{definition}

For the quantitative comparison estimate we temporarily impose
deterministic bounds on the viscous approximations.

\medskip
\noindent
\textbf{(A5$'$) Uniform bounds for the viscous families.}
In addition to \textbf{(A5)}, suppose that the corresponding viscous
families satisfy
\begin{equation}\label{eq:A5prime}
\operatorname*{ess\,sup}_{\omega}
\sup_{\varepsilon>0}
\sup_{t\in[0,T]}
\|u^\varepsilon(t)\|_{L^1(\R)}
\le
M_1,
\qquad
\operatorname*{ess\,sup}_{\omega}
\sup_{\varepsilon>0}
\sup_{t\in[0,T]}
\TV(u^\varepsilon(t))
\le
M_{\mathrm{BV}},
\end{equation}
and analogously for the approximations $v^{\varepsilon'}$.

\begin{theorem}[$L^1$-stability and uniqueness]
\label{thm:uniq}
Assume \textbf{(A1)}--\textbf{(A5)},
\textbf{(K1)}--\textbf{(K3)}, and \textbf{(A5$'$)}.
Let $u$ and $v$ be vanishing-viscosity solutions of
\eqref{eq:u}--\eqref{eq:id}, defined on the same stochastic basis and
driven by the same Brownian motion and Poisson random measure, with
initial data $u_0$ and $v_0$, respectively. Then, for a.e.\
$t\in[0,T]$,
\begin{equation}\label{eq:L1-stability-main}
\mathbb E
\|u(t)-v(t)\|_{L^1(\R)}
\le
e^{\mathcal Ct}
\mathbb E
\|u_0-v_0\|_{L^1(\R)},
\end{equation}
where one may take
\begin{equation}\label{eq:stability-constant}
\mathcal C
=
2\widetilde L_G
+
2L_fM_{\mathrm{BV}}\|\mu\|_{L^\infty}
+
3L_fM_1\|\mu'\|_{L^\infty}.
\end{equation}
In particular, if $u_0=v_0$ a.s., then
\[
u=v
\qquad
\text{a.s. on }
\Omega\times(0,T)\times\R .
\]
\end{theorem}



\begin{remark}[Why the comparison is performed at equal times]
\label{rem:equal-time}
The equal-time formulation is essential for both stochastic
contributions. If entropy inequalities for $u$ and $v$ were derived
separately and subsequently added, the Brownian corrections would
involve
\[
\Phi(u)^2+\Phi(v)^2,
\]
rather than the useful difference
\[
\big(\Phi(u)-\Phi(v)\big)^2.
\]
Similarly, the two jump increments would remain separated.

In contrast, the mixed second derivative of
\[
(p,q)
\longmapsto
\beta_\delta(p-q)
\]
and the simultaneous jumps of the pair $(u,v)$ produce directly
\[
\big(\Phi(u)-\Phi(v)\big)^2
\qquad\text{and}\qquad
G(u,z)-G(v,z).
\]
These structures allow the Lipschitz assumptions on $\Phi$ and $G$ to
be used directly in terms of $|u-v|$. This is the stochastic analogue
of the equal-time comparison mechanism used in
\cite{FengNualart2008,BVW2012,BKM2015}.
\end{remark}

\begin{proof}[Proof of Theorem~\ref{thm:uniq}]

We compare the viscous approximations
$u^\varepsilon$ and $v^{\varepsilon'}$ and obtain estimates independent
of the regularization parameters. The proof is divided into five
steps.

\medskip
\noindent
\textbf{Step 1: regularized entropy and doubling kernel.}

For $\delta>0$, let
\begin{equation}\label{eq:beta-delta-props}
\beta_\delta(r)
:=
\sqrt{r^2+\delta^2}-\delta .
\end{equation}
Then
\[
\beta_\delta'(r)
=
\frac{r}{\sqrt{r^2+\delta^2}},
\qquad
\beta_\delta''(r)
=
\frac{\delta^2}
{(r^2+\delta^2)^{3/2}}.
\]
In particular, $\beta_\delta$ is smooth, convex and even,
$\beta_\delta(0)=0$, and
\begin{equation}\label{eq:beta-props-list}
|r|-\delta
\le
\beta_\delta(r)
\le
|r|,
\qquad
|\beta_\delta'(r)|
\le1,
\qquad
\beta_\delta'(0)=0,
\qquad
\int_\R\beta_\delta''(r)\,\d r=2 .
\end{equation}
Moreover,
\begin{equation}\label{eq:beta-remainder}
0
\le
\beta_\delta''(r)r^2
\le
\min
\left\{
\frac{2}{3\sqrt3}\delta,\,
|r|
\right\},
\qquad
\big|
\beta_\delta(r)-r\beta_\delta'(r)
\big|
\le
\delta .
\end{equation}
The normalization $\beta_\delta(0)=0$ is important on the unbounded
domain $\R$: unlike the more common approximation
$\sqrt{r^2+\delta^2}$, it satisfies
\[
0\le\beta_\delta(r)\le|r|,
\]
and therefore
$\beta_\delta(w)\in L^1(\R)$ whenever $w\in L^1(\R)$.

Let $\rho\in C_c^\infty(\R)$ be nonnegative and even, with
\[
\int_\R\rho(r)\,\d r=1,
\qquad
\supp\rho\subset[-1,1],
\]
and define
\[
\rho_\kappa(r)
=
\frac1\kappa
\rho\!\left(\frac r\kappa\right),
\qquad
C_\rho
:=
\int_\R|\rho'(r)|\,\d r .
\]
Then
\[
\int_\R|\rho_\kappa'(r)|\,\d r
=
\frac{C_\rho}{\kappa}.
\]

We also introduce a spatial cut-off. Let
$\chi\in C_c^\infty(\R)$ satisfy
\[
0\le\chi\le1,
\qquad
\chi=1
\quad\text{on }[-1,1],
\]
and put $\chi_R(x)
= \chi(x/R).$ The doubled test function is
\begin{equation}\label{eq:phi-eps-def}
\phi(x,y)
:=
\chi_R(x)\rho_\kappa(x-y).
\end{equation}
Its derivatives satisfy
\[
\partial_x\phi
=
\chi_R'(x)\rho_\kappa(x-y)
+
\chi_R(x)\rho_\kappa'(x-y),
\]
and
\[
\partial_y\phi
=
-\chi_R(x)\rho_\kappa'(x-y).
\]
Consequently,
\begin{equation}\label{eq:phi-diagonal}
\partial_x\phi+\partial_y\phi
=
\chi_R'(x)\rho_\kappa(x-y).
\end{equation}
Identity \eqref{eq:phi-diagonal} is the basic cancellation used in the
transport estimate.

Let now $u^\varepsilon$ and $v^{\varepsilon'}$ solve the viscous
problems with viscosities $\varepsilon,\varepsilon'>0$. For notational
simplicity, write
\[
u=u^\varepsilon(t,x),
\qquad
v=v^{\varepsilon'}(t,y),
\]
and
\[
\mathcal U
=
\mathcal U[u^\varepsilon](t,x),
\qquad
\mathcal V
=
\mathcal U[v^{\varepsilon'}](t,y).
\]

From Section~\ref{subsec:kernel} and \eqref{eq:A5prime},
\begin{align}
\|\mathcal U[u^\varepsilon]
-\mathcal U[v^{\varepsilon'}]\|_{L^\infty}
&\le
\|\mu\|_{L^\infty}
\|u^\varepsilon-v^{\varepsilon'}\|_{L^1},
\label{eq:U-stab-used}
\\
\|\partial_x\mathcal U[u^\varepsilon]
-\partial_x\mathcal U[v^{\varepsilon'}]\|_{L^\infty}
&\le
\|\mu'\|_{L^\infty}
\|u^\varepsilon-v^{\varepsilon'}\|_{L^1},
\nonumber\\
\|\mathcal U[u^\varepsilon]\|_{L^\infty}
&\le
\|\mu\|_{L^\infty}M_1,
\nonumber\\
\|\partial_x\mathcal U[u^\varepsilon]\|_{L^\infty}
&\le
\|\mu'\|_{L^\infty}M_1,
\label{eq:U-uniform-used}
\end{align}
and the same bounds hold for $v^{\varepsilon'}$.

\medskip
\noindent
\textbf{Step 2: equal-time It\^o--L\'evy formula.}

Define
\[
X(t)
:=
\int_{\R^2}
\beta_\delta
\big(
u^\varepsilon(t,x)-v^{\varepsilon'}(t,y)
\big)
\phi(x,y)\,\d x\,\d y.
\]
Applying the It\^o--L\'evy formula to the pair
$(u^\varepsilon(t,x),v^{\varepsilon'}(t,y))$ gives
\begin{equation}\label{eq:ito-pair}
X(t)
=
X(0)
+
\int_0^t
\big(
\mathcal T+\mathcal E+\mathcal B+\mathcal J
\big)\,\d s
+
\mathcal M^W(t)
+
\mathcal M^\eta(t),
\end{equation}
where
\begin{align}
\mathcal T
&:=
\int_{\R^2}
\Big[
-\beta_\delta'(u-v)
\partial_x(f(u)\mathcal U)
+
\beta_\delta'(u-v)
\partial_y(f(v)\mathcal V)
\Big]
\phi\,\d x\,\d y,
\label{eq:transport-def}
\\
\mathcal E
&:=
\int_{\R^2}
\beta_\delta'(u-v)
\Big(
\varepsilon\partial_{xx}u
-
\varepsilon'\partial_{yy}v
\Big)
\phi\,\d x\,\d y,
\label{eq:visc-def}
\\
\mathcal B
&:=
\frac12
\int_{\R^2}
\beta_\delta''(u-v)
\big(
\Phi(u)-\Phi(v)
\big)^2
\phi\,\d x\,\d y,
\label{eq:brown-def}
\\
\mathcal J
&:=
\int_{\R^2}\!\!\int_Z
\Big[
\beta_\delta(u-v+\Delta G)
-
\beta_\delta(u-v)
-
\beta_\delta'(u-v)\Delta G
\Big]
\nu(\d z)\,\d x\,\d y,
\label{eq:jump-def}
\end{align}
with
\[
\Delta G
:=
G(u,z)-G(v,z).
\]

Notice that the Brownian correction in \eqref{eq:brown-def} contains
the square of the \emph{difference} of the two noise coefficients.
Likewise, the jump correction in \eqref{eq:jump-def} involves the
single increment $\Delta G$. This is precisely the advantage of
applying the It\^o--L\'evy formula to the pair rather than to the two
solutions separately.

The corresponding martingales are
\[
\mathcal M^W(t)
=
\int_0^t
\int_{\R^2}
\beta_\delta'(u-v)
\big(
\Phi(u)-\Phi(v)
\big)
\phi\,\d x\,\d y\,\d W(s)
\]
and
\begin{align*}
\mathcal M^\eta(t)
&=
\int_0^t\!\!\int_Z
\int_{\R^2}
\Big[
\beta_\delta(u_- -v_-+\Delta G_-)
-
\beta_\delta(u_- -v_-)
\Big] \times\phi(x,y)\,\d x\,\d y\,
\widetilde\eta(\d z,\d s).
\end{align*}
The compact spatial cut-off, \eqref{eq:A5prime}, and assumptions
\textbf{(A2)}--\textbf{(A3)} guarantee their integrability. Therefore
\[
\mathbb E\mathcal M^W(t)
=
\mathbb E\mathcal M^\eta(t)
=
0,
\]
and taking expectations in \eqref{eq:ito-pair} gives
\begin{equation}\label{eq:main-after-expect}
\mathbb E X(t)
=
\mathbb E X(0)
+
\int_0^t
\mathbb E
\big(
\mathcal T+\mathcal E+\mathcal B+\mathcal J
\big)\,\d s .
\end{equation}

\medskip
\noindent
\textbf{Step 3: viscous, Brownian, and jump contributions.}

We first estimate the terms that do not involve the nonlocal transport
structure.

Integrating \eqref{eq:visc-def} by parts in $x$ and $y$ yields
\begin{align}
\mathcal E
&=
-\int_{\R^2}
\beta_\delta''(u-v)
\Big(
\varepsilon|\partial_xu|^2
+
\varepsilon'|\partial_yv|^2
\Big)
\phi\,\d x\,\d y -\int_{\R^2}
\beta_\delta'(u-v)
\Big(
\varepsilon\partial_xu\,\partial_x\phi
-
\varepsilon'\partial_yv\,\partial_y\phi
\Big)
\,\d x\,\d y .
\label{eq:visc-split}
\end{align}
The first integral in \eqref{eq:visc-split} is nonpositive by the
convexity of $\beta_\delta$ and may therefore be discarded in an upper
estimate. Using
\[
|\beta_\delta'|\le1,
\qquad
\TV(u^\varepsilon),
\TV(v^{\varepsilon'})
\le M_{\mathrm{BV}},
\]
together with
\[
\|\chi_R'\|_\infty
=
R^{-1}\|\chi'\|_\infty,
\qquad
\|\rho_\kappa'\|_{L^1}
=
C_\rho/\kappa,
\]
we obtain
\begin{equation}\label{eq:visc-err}
\mathcal E
\le
M_{\mathrm{BV}}
\left(
\frac{\varepsilon\|\chi'\|_{L^\infty}}{R}
+
\frac{C_\rho(\varepsilon+\varepsilon')}{\kappa}
\right).
\end{equation}

For the Brownian correction, the Lipschitz continuity of $\Phi$ gives
\[
|\Phi(u)-\Phi(v)|
\le
L_\Phi|u-v|.
\]
Hence
\[
\mathcal B
\le
\frac{L_\Phi^2}{2}
\int_{\R^2}
\beta_\delta''(u-v)
|u-v|^2
\phi\,\d x\,\d y.
\]
By \eqref{eq:beta-remainder}, $\beta_\delta''(r)r^2
\le C\delta,$ and therefore, for fixed $R$ and $\kappa$,
\begin{equation}\label{eq:brown-vanish}
\lim_{\delta\downarrow0}
\mathbb E
\int_0^t
\mathcal B(s)\,\d s
=
0.
\end{equation}
This is stated separately as Lemma~\ref{lem:BM-good}.

For the jump compensator, convexity gives
\[
0
\le
\beta_\delta(a+b)
-
\beta_\delta(a)
-
\beta_\delta'(a)b,
\]
while $|\beta_\delta'|\le1$ implies
\[
\beta_\delta(a+b)
-
\beta_\delta(a)
-
\beta_\delta'(a)b
\le
2|b|.
\]
Taking
\[
a=u-v,
\qquad
b=\Delta G,
\]
and using the $L^1(Z,\nu)$-Lipschitz assumption in
\textbf{(A3)}, we obtain
\begin{equation}\label{eq:jump-bound}
0
\le
\mathcal J
\le
2\widetilde L_G
\int_{\R^2}
|u-v|
\phi\,\d x\,\d y .
\end{equation}
Unlike the Brownian correction, the jump compensator does not in
general vanish as $\delta\downarrow0$; instead, it contributes to the
Gr\"onwall coefficient.

\medskip
\noindent
\textbf{Step 4: the nonlocal transport contribution.}

This is the main part of the proof. The shifted entropy flux
\[
\Psi_k(r)
=
\int_k^r
\beta_\delta'(s-k)f'(s)\,\d s
\]
is designed precisely to handle the spatial dependence of the
nonlocal velocity.

Integrating the two terms in \eqref{eq:transport-def} by parts and
using the entropy-flux identity of Section~\ref{sec:entropy}, we obtain
\[
\mathcal T
=
\mathfrak F-\mathcal S,
\]
where
\begin{align}
\mathfrak F
&:=
\int_{\R^2}
\Big[
\Psi_v(u)\mathcal U\,\partial_x\phi
+
\Psi_u(v)\mathcal V\,\partial_y\phi
\Big]
\,\d x\,\d y,
\label{eq:F-def}
\\
\mathcal S
&:=
\int_{\R^2}
\Big(
[\beta_\delta'(u-v)f(u)-\Psi_v(u)]
\partial_x\mathcal U +
[\beta_\delta'(v-u)f(v)-\Psi_u(v)]
\partial_y\mathcal V
\Big)
\phi\,\d x\,\d y .
\label{eq:S-def}
\end{align}

We treat these two terms separately.

\medskip
\noindent
\emph{Step 4a: decomposition of $\mathfrak F$.}

From \eqref{eq:phi-diagonal}, $\partial_x\phi
= \chi_R'\rho_\kappa -
\partial_y\phi.$ Therefore,
\begin{equation}\label{eq:F-split}
\mathfrak F
=
\mathfrak F_1+\mathfrak F_2+\mathfrak F_3,
\end{equation}
where
\[
\mathfrak F_1
:=
\int_{\R^2}
\Psi_v(u)\mathcal U
\chi_R'(x)\rho_\kappa(x-y)
\,\d x\,\d y, \qquad \mathfrak F_2
:=
\int_{\R^2}
\Psi_v(u)
(\mathcal V-\mathcal U)
\partial_y\phi
\,\d x\,\d y, 
\]
and
\[
\mathfrak F_3
:=
\int_{\R^2}
\big[
\Psi_u(v)-\Psi_v(u)
\big]
\mathcal V
\partial_y\phi
\,\d x\,\d y.
\]

The first term is a pure boundary contribution. Since
\[
|\Psi_v(u)|
\le
L_f|u-v|, \quad \|\mathcal U\|_\infty
\le
M_1\|\mu\|_\infty, \quad \mbox{and} \quad \|\chi_R'\|_\infty
=
\frac{\|\chi'\|_\infty}{R},
\]
we obtain
\begin{equation}\label{eq:F1-bound}
|\mathfrak F_1|
\le
\frac{
2L_fM_1^2
\|\mu\|_{L^\infty}
\|\chi'\|_{L^\infty}
}{R}.
\end{equation}
Hence
\[
\mathfrak F_1
\longrightarrow0
\qquad\text{as }R\uparrow\infty.
\]

\medskip
\noindent
\emph{Step 4b: estimate of $\mathfrak F_2$.}

The derivative $\partial_y\phi$ contains the singular factor
$\rho_\kappa'$, whose $L^1$-norm is of order $\kappa^{-1}$.
Estimating $\mathfrak F_2$ directly would therefore lose a factor
$\kappa^{-1}$. To avoid this, we integrate by parts with respect to
$y$:
\[
\mathfrak F_2
=
-\int_{\R^2}
\Big\{
\partial_y[\Psi_v(u)]
(\mathcal V-\mathcal U)
+
\Psi_v(u)\partial_y\mathcal V
\Big\}
\phi\,\d x\,\d y.
\]

Since $u=u(t,x)$ is independent of $y$, the $y$-dependence of
$\Psi_v(u)$ comes only through $v=v(t,y)$. Differentiating with
respect to the lower entropy level gives
\[
\partial_v\Psi_v(u)
=
-\beta_\delta'(0)f'(v)
-
\int_v^u
\beta_\delta''(r-v)f'(r)\,\d r .
\]
Because $\beta_\delta'(0)=0$,
\[
\partial_v\Psi_v(u)
=
-\int_v^u
\beta_\delta''(r-v)f'(r)\,\d r .
\]
Consequently,
\[
|\partial_v\Psi_v(u)|
\le
L_f
\int_\R
\beta_\delta''(s)\,\d s
=
2L_f,
\]
and therefore
\begin{equation}\label{eq:Psi-v-derivative}
|\partial_y[\Psi_v(u)]|
\le
2L_f|\partial_yv|.
\end{equation}

We next estimate $\mathcal V-\mathcal U$ on the support of
$\rho_\kappa$. Write
\[
\mathcal V(t,y)-\mathcal U(t,x)
=
\big[
\mathcal U[v](t,y)-\mathcal U[v](t,x)
\big]
+
\big[
\mathcal U[v](t,x)-\mathcal U[u](t,x)
\big].
\]
Since $|x-y|\le\kappa$ on
$\supp\rho_\kappa(x-y)$,
\[
|
\mathcal U[v](t,y)-\mathcal U[v](t,x)
|
\le
|x-y|
\|\partial_x\mathcal U[v]\|_\infty
\le
\kappa
M_1\|\mu'\|_\infty.
\]
On the other hand,
\[
\|
\mathcal U[v]-\mathcal U[u]
\|_\infty
\le
\|\mu\|_\infty
\|u^\varepsilon-v^{\varepsilon'}\|_{L^1}.
\]
Hence
\begin{equation}\label{eq:UV-pointwise}
|\mathcal V-\mathcal U|
\le
\|\mu\|_{L^\infty}
\|u^\varepsilon-v^{\varepsilon'}\|_{L^1}
+
\kappa
M_1\|\mu'\|_{L^\infty}.
\end{equation}

Using \eqref{eq:Psi-v-derivative},
\eqref{eq:UV-pointwise},
$\TV(v^{\varepsilon'})\le M_{\mathrm{BV}}$, and
\[
|\Psi_v(u)|
\le
L_f|u-v|,
\qquad
\|\partial_y\mathcal V\|_\infty
\le
M_1\|\mu'\|_\infty,
\]
we arrive at
\begin{align}
|\mathfrak F_2|
&\le
2L_fM_{\mathrm{BV}}
\left[
\|\mu\|_{L^\infty}
\|u^\varepsilon-v^{\varepsilon'}\|_{L^1}
+
\|\mu'\|_{L^\infty}M_1\kappa
\right] +
L_f
\|\mu'\|_{L^\infty}M_1
\int_{\R^2}
|u-v|\phi\,\d x\,\d y .
\label{eq:mixed-final}
\end{align}
Thus the apparently singular derivative of $\rho_\kappa$ has been
completely removed; only an $O(\kappa)$ error remains.

\medskip
\noindent
\emph{Step 4c: estimate of $\mathfrak F_3$.}

The term $\mathfrak F_3$ measures the lack of exact antisymmetry of the
regularized shifted entropy flux. This term disappears in the
Kru\v{z}kov limit $\delta\downarrow0$.

Set
\[
h_\delta(r)
:=
\frac{r}{\sqrt{r^2+\delta^2}},
\qquad r\ge0.
\]
Suppose first that $u\ge v$ and put $d=u-v$. Then
\[
\Psi_v(u)
=
\int_0^d
h_\delta(s)
f'(v+s)\,\d s.
\]
Changing variables in $\Psi_u(v)$ gives
\[
\Psi_u(v)
=
\int_0^d
h_\delta(d-s)
f'(v+s)\,\d s.
\]
Hence
\begin{align*}
|\Psi_u(v)-\Psi_v(u)|
&\le
L_f
\int_0^d
|h_\delta(d-s)-h_\delta(s)|\,\d s \le
2L_f
\int_0^\infty
\big(
1-h_\delta(s)
\big)\,\d s.
\end{align*}
Since
\[
\int_0^\infty
\left(
1-
\frac{s}{\sqrt{s^2+\delta^2}}
\right)\d s
=
\delta,
\]
we obtain
\begin{equation}\label{eq:Psi-symmetry}
|\Psi_u(v)-\Psi_v(u)|
\le
2L_f\delta .
\end{equation}
The same estimate holds when $u<v$.

Using $\|\mathcal V\|_\infty
\le M_1\|\mu\|_\infty,$ together with
\[
\int_{\R^2}
|\partial_y\phi|\,\d x\,\d y
\le
C\frac{R}{\kappa},
\]
we obtain
\begin{equation}\label{eq:F3-bound}
|\mathfrak F_3|
\le
C
L_f
\|\mu\|_{L^\infty}
M_1
\frac{R}{\kappa}\delta .
\end{equation}
Thus, for fixed $R$ and $\kappa$,
\[
\mathfrak F_3
\longrightarrow0
\qquad
\text{as }\delta\downarrow0.
\]

\medskip
\noindent
\emph{Step 4d: the correction term $\mathcal S$.}

We now examine \eqref{eq:S-def}. Since $\beta_\delta'(r)
\longrightarrow
\sgn(r)$ and
\[
\Psi_v(u)
\longrightarrow
\sgn(u-v)
\big(
f(u)-f(v)
\big),
\]
we have
\[
\beta_\delta'(u-v)f(u)-\Psi_v(u)
\longrightarrow
\sgn(u-v)f(v).
\]
Similarly,
\[
\beta_\delta'(v-u)f(v)-\Psi_u(v)
\longrightarrow
-\sgn(u-v)f(u).
\]
Therefore, after first letting $\delta\downarrow0$ and then identifying
the diagonal $x=y$ by $\kappa\downarrow0$,
\begin{equation}\label{eq:limit-S}
\mathcal S
\longrightarrow
\int_\R
\sgn(u-v)
\Big[
f(v)\partial_x\mathcal U[u]
-
f(u)\partial_x\mathcal U[v]
\Big]
\chi_R\,\d x .
\end{equation}

To estimate this term, write
\begin{align*}
f(v)\partial_x\mathcal U[u]
-
f(u)\partial_x\mathcal U[v]
&=
f(v)
\Big(
\partial_x\mathcal U[u]
- \partial_x\mathcal U[v] \Big) +
\big( f(v)-f(u)
\big)
\partial_x\mathcal U[v].
\end{align*}
For the first contribution,
\[
\int_\R
|f(v)|
\,
|
\partial_x\mathcal U[u]
-
\partial_x\mathcal U[v]
|
\,\d x
\]
is bounded by
\[
L_f\|v\|_{L^1}
\|\mu'\|_\infty
\|u-v\|_{L^1}
\le
L_fM_1
\|\mu'\|_\infty
\|u-v\|_{L^1}.
\]
For the second,
\[
\int_\R
|f(u)-f(v)|
|\partial_x\mathcal U[v]|
\,\d x
\le
L_fM_1
\|\mu'\|_\infty
\|u-v\|_{L^1}.
\]
Consequently,
\begin{equation}\label{eq:S-est2}
\left|
\lim_{\kappa\downarrow0}
\lim_{\delta\downarrow0}
\mathcal S
\right|
\le
2L_fM_1
\|\mu'\|_{L^\infty}
\|u^\varepsilon-v^{\varepsilon'}\|_{L^1(\R)} .
\end{equation}

\medskip
\noindent
\textbf{Step 5: passage to the limits and Gr\"onwall estimate.}

We now combine the estimates above. The order of the limits is
important:
\[
\delta\downarrow0,
\qquad
\varepsilon,\varepsilon'\downarrow0,
\qquad
\kappa\downarrow0,
\qquad
R\uparrow\infty.
\]

First fix $\varepsilon,\varepsilon',\kappa$ and $R$ and let
$\delta\downarrow0$. Then \eqref{eq:brown-vanish} removes the Brownian
correction, while \eqref{eq:F3-bound} removes $\mathfrak F_3$.

Next, with $\kappa$ and $R$ fixed, let
$\varepsilon,\varepsilon'\downarrow0$. The viscosity estimate
\eqref{eq:visc-err} then vanishes. The convergence in
Definition~\ref{def:vv}, together with the Lipschitz continuity of the
coefficients and the convolution estimates, permits passage to the
remaining terms on the support of $\chi_R$.

We then send $\kappa\downarrow0$. Since
$\rho_\kappa$ is an approximate identity,
translation continuity in $L^1(\R)$ gives
\begin{align}
&
\int_{\R^2}
|u(t,x)-v(t,y)|
\rho_\kappa(x-y)
\chi_R(x)\,\d x\,\d y \longrightarrow
\int_\R
|u(t,x)-v(t,x)|
\chi_R(x)\,\d x .
\label{eq:diag-limit}
\end{align}
At the same time, the explicit $O(\kappa)$ term in
\eqref{eq:mixed-final} disappears.

Finally, we let $R\uparrow\infty$. Since
$0\le\chi_R\le1$ and $\chi_R(x)\to1$ pointwise,
monotone or dominated convergence yields the full $L^1(\R)$ norm,
while the boundary term $\mathfrak F_1$ vanishes by
\eqref{eq:F1-bound}.

Collecting the jump contribution \eqref{eq:jump-bound}, the estimate
for $\mathfrak F_2$, and \eqref{eq:S-est2}, we obtain, for a.e.\
$t\in[0,T]$,
\begin{equation}\label{eq:L1-ineq}
\mathbb E
\|u(t)-v(t)\|_{L^1(\R)}
\le
\mathbb E
\|u_0-v_0\|_{L^1(\R)}
+
\mathcal C
\int_0^t
\mathbb E
\|u(s)-v(s)\|_{L^1(\R)}
\,\d s,
\end{equation}
where
\[
\mathcal C
=
2\widetilde L_G
+
2L_fM_{\mathrm{BV}}\|\mu\|_{L^\infty}
+
3L_fM_1\|\mu'\|_{L^\infty}.
\]

Indeed, $2\widetilde L_G$ comes from the jump compensator, and $2L_fM_{\mathrm{BV}}\|\mu\|_\infty$ comes from the first term in \eqref{eq:mixed-final}, while the three
contributions involving
$L_fM_1\|\mu'\|_\infty$ arise from the second part of
\eqref{eq:mixed-final} and the two terms in \eqref{eq:S-est2}.

Gronwall's inequality now yields
\[
\mathbb E
\|u(t)-v(t)\|_{L^1(\R)}
\le
e^{\mathcal Ct}
\mathbb E
\|u_0-v_0\|_{L^1(\R)}.
\]
This proves \eqref{eq:L1-stability-main}.
\end{proof}

\subsection{A stopped Cauchy estimate}

The preceding proof also yields a quantitative estimate between two
viscous approximations. This will be used later to remove the viscosity.

\begin{corollary}[Stopped Cauchy estimate]
\label{cor:cauchy}
Let $u^\varepsilon$ and $u^{\varepsilon'}$ be viscous approximations
corresponding to the same limiting initial datum, with regularized
initial values $u_0^\varepsilon$ and $u_0^{\varepsilon'}$.
Let $\sigma$ be a stopping time such that on $[0,\sigma]$,
\[
\|u^\varepsilon(t)\|_{L^1}
+
\|u^{\varepsilon'}(t)\|_{L^1}
\le M_1,
\]
and
\[
\TV(u^\varepsilon(t))
+
\TV(u^{\varepsilon'}(t))
\le M_{\mathrm{BV}}
\qquad\text{a.s.}
\]
Then
\begin{align}
\sup_{0\le t\le T}
\mathbb E
\left\|
u^\varepsilon(t\wedge\sigma)
-
u^{\varepsilon'}(t\wedge\sigma)
\right\|_{L^1(\R)}
&\le
Ce^{\mathcal CT}
\Big[
\mathbb E
\|u_0^\varepsilon-u_0^{\varepsilon'}\|_{L^1} +
M_{\mathrm{BV}}
\sqrt{\varepsilon+\varepsilon'}
\Big].
\label{eq:cauchy-rate}
\end{align}
In particular, if
\[
\mathbb E
\|u_0^\varepsilon-u_0^{\varepsilon'}\|_{L^1}
=
O(\varepsilon+\varepsilon'),
\]
then
\[
\sup_{0\le t\le T}
\mathbb E
\left\|
u^\varepsilon(t\wedge\sigma)
-
u^{\varepsilon'}(t\wedge\sigma)
\right\|_{L^1}
=
O
\big(
\sqrt{\varepsilon+\varepsilon'}
\big).
\]
\end{corollary}

\begin{proof}
We repeat the proof of Theorem~\ref{thm:uniq} up to the stopping time
$\sigma$. Since all $L^1$- and $BV$-bounds are deterministic on this
interval, the same estimates apply.

After letting $\delta\downarrow0$, we first send
$R\uparrow\infty$, keeping $\kappa$, $\varepsilon$, and
$\varepsilon'$ fixed. This eliminates the boundary term without
requiring any spatial moment or tail assumption.

The remaining viscosity contribution satisfies
\[
\mathcal E
\le
C
M_{\mathrm{BV}}
\frac{\varepsilon+\varepsilon'}{\kappa}.
\]
On the other hand, the error caused by identifying $x$ and $y$ is of
order $\kappa$. Indeed, for $w\in BV(\R)$,
\begin{equation}\label{eq:BV-translation}
\|w(\cdot+h)-w\|_{L^1(\R)}
\le
|h|\TV(w).
\end{equation}
Since $\rho_\kappa$ is supported in $|x-y|\le\kappa$,
\eqref{eq:BV-translation} yields an error bounded by $C\kappa M_{\mathrm{BV}}.$

The initial doubled term also has to be treated carefully. It does not
vanish merely because the two viscous equations approximate the same
initial datum, since the doubled integral compares the values at
different points $x$ and $y$. We write
\begin{align*}
&
\int_{\R^2}
|u_0^\varepsilon(x)-u_0^{\varepsilon'}(y)|
\rho_\kappa(x-y)\,\d x\,\d y
\\
&\quad
\le
\int_{\R^2}
|u_0^\varepsilon(x)-u_0^{\varepsilon'}(x)|
\rho_\kappa(x-y)\,\d x\,\d y +
\int_{\R^2}
|u_0^{\varepsilon'}(x)-u_0^{\varepsilon'}(y)|
\rho_\kappa(x-y)\,\d x\,\d y
\\
&\quad
\le
\|u_0^\varepsilon-u_0^{\varepsilon'}\|_{L^1}
+
C\kappa
\TV(u_0^{\varepsilon'}).
\end{align*}

Thus the comparison estimate becomes
\begin{align*}
&
\mathbb E
\left\|
u^\varepsilon(t\wedge\sigma)
-
u^{\varepsilon'}(t\wedge\sigma)
\right\|_{L^1} \le
Ce^{\mathcal Ct}
\left[
\mathbb E
\|u_0^\varepsilon-u_0^{\varepsilon'}\|_{L^1}
+
CM_{\mathrm{BV}}
\left(
\frac{\varepsilon+\varepsilon'}{\kappa}
+
\kappa
\right)
\right].
\end{align*}
Optimizing with respect to $\kappa$ gives $\kappa
= \sqrt{\varepsilon+\varepsilon'},$ and hence \eqref{eq:cauchy-rate}.
\end{proof}

\smallskip

\begin{remark}[Mollified initial data]
If $u_0^\varepsilon = \rho_\varepsilon*u_0$ and $u_0\in BV(\R)$, then
\[
\|u_0^\varepsilon-u_0\|_{L^1}
\le
C\varepsilon\,\TV(u_0).
\]
Consequently,
\[
\|u_0^\varepsilon-u_0^{\varepsilon'}\|_{L^1}
\le
C(\varepsilon+\varepsilon')
\TV(u_0),
\]
so the assumption on the initial regularizations in
Corollary~\ref{cor:cauchy} is satisfied.
\end{remark}

\subsection{Removal of the deterministic bounds}

We now remove assumption \textbf{(A5$'$)} when the two solutions have
the same initial datum.

\begin{corollary}[Uniqueness without deterministic a priori bounds]
\label{cor:uniq-localised}
Assume
\textbf{(A1)}--\textbf{(A5)} and
\textbf{(K1)}--\textbf{(K4)}, but not
\textbf{(A5$'$)}. Let $u$ and $v$ be vanishing-viscosity solutions
with
\[
u_0=v_0
\qquad\text{a.s.}
\]
Then
\[
u=v
\qquad
\text{a.s. on }
\Omega\times(0,T)\times\R .
\]
\end{corollary}

\begin{proof}
Let $u^\varepsilon$ and $v^{\varepsilon'}$ be viscous approximations
converging respectively to $u$ and $v$. By
Definition~\ref{def:vv}, their initial regularizations are uniformly
$BV$-stable, so Proposition~\ref{prop:BV} applies to both viscous families.
For $M\ge1$, introduce
\[
\sigma_M^{\varepsilon,\varepsilon'}
:=
\inf
\Big\{
t\ge0:
\|u^\varepsilon(t)\|_{L^1}
+
\|v^{\varepsilon'}(t)\|_{L^1}
+
\TV(u^\varepsilon(t))
+
\TV(v^{\varepsilon'}(t))
\ge M
\Big\}
\wedge T .
\]
On
$[0,\sigma_M^{\varepsilon,\varepsilon'}]$,
the deterministic estimates required in the comparison argument hold
with constants depending only on $M$. Hence
Corollary~\ref{cor:cauchy} gives, for every fixed $M$,
\begin{align}
&
u^\varepsilon
\big(
\cdot\wedge
\sigma_M^{\varepsilon,\varepsilon'}
\big)
-
v^{\varepsilon'}
\big(
\cdot\wedge
\sigma_M^{\varepsilon,\varepsilon'}
\big) \longrightarrow0
\quad\text{in probability}
\label{eq:stopped-uniq-prob}
\end{align}
as
$\varepsilon,\varepsilon'\downarrow0$.

On the other hand, the uniform $L^1$-moment estimate and the localized
$BV$-estimate obtained in Section~\ref{sec:existence} imply
\begin{equation}\label{eq:uniq-exit-prob}
\lim_{M\uparrow\infty}
\sup_{\varepsilon,\varepsilon'>0}
\mathbb P
\big(
\sigma_M^{\varepsilon,\varepsilon'}<T
\big)
=
0.
\end{equation}
Indeed, the $L^1$ norms are uniformly bounded in probability by the
second-moment estimate, while
Proposition~\ref{prop:BV} gives the corresponding boundedness in
probability of the total variation.

Let $\lambda>0$. For any fixed $M$,
\begin{align*}
&
\mathbb P
\left(
\|u^\varepsilon(t)-v^{\varepsilon'}(t)\|_{L^1}
>\lambda
\right) \le
\mathbb P
\big(
\sigma_M^{\varepsilon,\varepsilon'}<T
\big)
+
\mathbb P
\left(
\left\|
u^\varepsilon
(t\wedge\sigma_M^{\varepsilon,\varepsilon'})
-
v^{\varepsilon'}
(t\wedge\sigma_M^{\varepsilon,\varepsilon'})
\right\|_{L^1}
>\lambda
\right).
\end{align*}
First let
$\varepsilon,\varepsilon'\downarrow0$ and use
\eqref{eq:stopped-uniq-prob}; then let $M\uparrow\infty$ and use
\eqref{eq:uniq-exit-prob}. We obtain
\[
\|u(t)-v(t)\|_{L^1}
=0
\qquad\text{a.s.}
\]
for a.e.\ $t$. Hence $u=v$ almost surely on
$\Omega\times(0,T)\times\R$.
\end{proof}

\begin{corollary}[Well-posedness in the vanishing-viscosity class]
\label{cor:wellposed}

Under the hypotheses of Theorem~\ref{thm:existence} and
Corollary~\ref{cor:uniq-localised}, problem
\eqref{eq:u}--\eqref{eq:id} admits a unique vanishing-viscosity
solution. By Section~\ref{subsec:eps-limit}, this solution is a
stochastic entropy solution in the sense of
Definition~\ref{solution-defn}.

If, in addition, the deterministic estimates
\textbf{(A5$'$)} hold, then the quantitative stability estimate
\eqref{eq:L1-stability-main} is valid for two solutions with possibly
different initial data.
\end{corollary}

The uniqueness established in this section is uniqueness within the
class of vanishing-viscosity solutions. It does not assert uniqueness
for every stochastic entropy solution satisfying
Definition~\ref{solution-defn}.

To extend the comparison principle to arbitrary entropy solutions, one
would have to compare an entropy solution directly with another
possibly nonsmooth entropy solution, or with a viscous approximation,
and handle the resulting random entropy level while preserving
adaptedness. Such an argument requires additional stochastic
doubling-of-variables machinery of the type developed in
\cite{BVW2012,BVW2014}, and is beyond the scope of the present work.

We also emphasize that the compactness involved in the Galerkin limit
$m\to\infty$ is logically independent of the comparison argument of
this section. The Galerkin limit is treated separately in
Section~\ref{subsec:galerkin-compactness}, whereas the Cauchy estimate
above is used only for the subsequent vanishing-viscosity limit
$\varepsilon\downarrow0$.

\medskip

\section{Auxiliary Estimates for the Stability Argument}
\label{sec:tech}\label{sec:lemmas}

The deterministic convolution estimates required above were already
established in Section~\ref{subsec:kernel}, and the spatial cut-off term
is controlled directly by \eqref{eq:F1-bound}. We therefore record only
the two stochastic estimates that are specific to the equal-time
It\^o--L\'evy comparison.

\subsection{Brownian correction}

\begin{lemma}[Vanishing of the Brownian correction]
\label{lem:BM-good}
Let $\beta_\delta$ be defined by
\eqref{eq:beta-delta-props}. Under \textbf{(A2)}, for
$\phi=\chi_R\rho_\kappa$ as in \eqref{eq:phi-eps-def},
\begin{equation}\label{eq:BM-good-lemma}
0
\le
\mathcal B
:=
\frac12
\int_{\R^2}
\beta_\delta''(u-v)
\big(\Phi(u)-\Phi(v)\big)^2
\phi\,\d x\,\d y
\le
C L_\Phi^2 R\delta ,
\end{equation}
where $C$ depends only on the fixed cut-off $\chi$.
Consequently,
\[
\mathcal B\longrightarrow0
\qquad\text{as }\delta\downarrow0
\]
for fixed $R$ and $\kappa$, uniformly with respect to the viscous
parameters.
\end{lemma}

\begin{proof}
By \textbf{(A2)},
\[
|\Phi(u)-\Phi(v)|
\le
L_\Phi|u-v|,
\]
and \eqref{eq:beta-remainder} gives
\[
\beta_\delta''(u-v)|u-v|^2
\le
\frac{2}{3\sqrt3}\delta .
\]
Therefore
\[
\mathcal B
\le
C L_\Phi^2\delta
\int_{\R^2}\phi(x,y)\,\d x\,\d y.
\]
Since $\int_\R\rho_\kappa(x-y)\,\d y=1$,
\[
\int_{\R^2}\phi(x,y)\,\d x\,\d y
=
\int_\R\chi_R(x)\,\d x
\le CR ,
\]
which proves the result.
\end{proof}

\subsection{Jump correction}

\begin{lemma}[Convexity estimate for the jump remainder]
\label{lem:jump-convex}
Let $\beta_\delta$ be convex with
$|\beta_\delta'|\le1$. Then, for every $a,b\in\R$,
\begin{equation}\label{eq:jump-two-sided}
0
\le
\beta_\delta(a+b)
-
\beta_\delta(a)
-
\beta_\delta'(a)b
\le
2|b|.
\end{equation}
Consequently, under \textbf{(A3)}, with
$\Delta G=G(u,z)-G(v,z)$,
\begin{equation}\label{eq:jump-good-lemma}
0
\le
\mathcal J
\le
2\widetilde L_G
\int_{\R^2}
|u-v|\phi\,\d x\,\d y .
\end{equation}
\end{lemma}

\begin{proof}
Convexity yields $\beta_\delta(a+b)
\ge
\beta_\delta(a)+\beta_\delta'(a)b.$ Moreover, $|\beta_\delta'|\le1$ implies
\[
|\beta_\delta(a+b)-\beta_\delta(a)|
\le|b|,
\qquad
|\beta_\delta'(a)b|
\le|b|,
\]
which proves \eqref{eq:jump-two-sided}. Taking
$a=u-v$ and $b=\Delta G$, integrating with respect to $\nu$, and using
\[
\int_Z
|G(u,z)-G(v,z)|\,\nu(\d z)
\le
\widetilde L_G|u-v|
\]
gives \eqref{eq:jump-good-lemma}.
\end{proof}

\begin{remark}
Unlike the Brownian correction, the jump compensator does not generally
vanish as $\delta\downarrow0$. Its contribution is instead bounded
linearly by the $L^1$-distance between the two solutions and is absorbed
in the final Gr\"onwall estimate. This is the role of the
$L^1(Z,\nu)$-Lipschitz condition in \textbf{(A3)}.
\end{remark}

\medskip

\section{Existence of Stochastic Entropy Solutions}
\label{sec:existence}

We now establish the existence of a stochastic entropy solution to
\eqref{eq:u}--\eqref{eq:id} by means of a vanishing-viscosity
approximation. The proof involves two successive limiting procedures.
For each fixed viscosity parameter $\varepsilon>0$, we first construct a
solution of the corresponding parabolic problem through a
Faedo--Galerkin approximation. We then derive estimates that are
uniform in $\varepsilon$ and pass to the limit $\varepsilon\downarrow0$.

There are two points that require particular care. First, the transport
operator
\[
u\longmapsto
\partial_x\big(f(u)\mathcal U[u]\big)
\]
is both nonlinear and nonlocal. Second, the natural $L^1$- and
$BV$-estimates cannot be justified directly at the Galerkin level,
because the nonlinear test functions used in these estimates do not
belong to the finite-dimensional Galerkin space and the projection does
not commute with nonlinear superposition. We therefore proceed in the
following order:
\[
\text{stopped Galerkin approximation}
\;\longrightarrow\;
\text{viscous solution}
\;\longrightarrow\;
L^1\text{- and }BV\text{-estimates}
\;\longrightarrow\;
\varepsilon\downarrow0 .
\]
This ordering avoids any circular dependence between the energy,
$L^1$, and $BV$ estimates.

\subsection{Viscous approximation}
\label{subsec:viscous}

We first regularize the initial datum.

\begin{lemma}[Approximation of the initial datum]
\label{lem:initial-approx}
There exists a family of $\mathcal F_0$-measurable random variables
$\{u_0^\varepsilon\}_{\varepsilon>0}$ such that
\[
u_0^\varepsilon
\in
L^2\!\left(
\Omega;
C^\infty(\R)\cap L^1(\R)\cap L^2(\R)
\right),
\]
and
\begin{align}
u_0^\varepsilon
&\longrightarrow
u_0
\quad\text{in }
L^2\!\left(
\Omega;L^1(\R)\cap L^2(\R)
\right),
\label{eq:initconv}
\\
\|u_0^\varepsilon\|_{L^1}
&\le
\|u_0\|_{L^1},
\label{eq:initL1}
\\
\|u_0^\varepsilon\|_{L^2}
&\le
\|u_0\|_{L^2}
\qquad\text{a.s.}
\label{eq:initL2}
\end{align}
If $u_0\in BV(\R)$ almost surely, then
\begin{equation}\label{eq:initBV}
\TV(u_0^\varepsilon)
\le
\TV(u_0)
\qquad\text{a.s.}
\end{equation}
A canonical choice is
\[
u_0^\varepsilon
=
\rho_\varepsilon*u_0,
\]
where $\rho_\varepsilon$ is a standard Friedrichs mollifier.
\end{lemma}

\begin{proof}
The convergence and $L^p$-stability properties follow from the standard
approximation properties of Friedrichs mollifiers. The variation
estimate follows from
\[
\TV(\rho_\varepsilon*u_0)
\le
\TV(u_0).
\]
\end{proof}

For $\varepsilon>0$, we consider the viscous problem
\begin{align}
\label{eq:viscous}
\d u^\varepsilon
+
\partial_x
\Big(
f(u^\varepsilon)\mathcal U[u^\varepsilon]
\Big)\,\d t
&=
\varepsilon\partial_{xx}u^\varepsilon\,\d t
+
\Phi(u^\varepsilon)\,\d W(t)
+
\int_Z
G(u^\varepsilon,z)\,
\widetilde\eta(\d z,\d t),
\\
u^\varepsilon(0)
&=
u_0^\varepsilon .
\label{eq:viscousIC}
\end{align}
The artificial viscosity is introduced only in the deterministic part
of the equation; the Brownian and jump coefficients, as well as the
nonlocal field
\[
\mathcal U[u^\varepsilon]
=
\mu*u^\varepsilon,
\]
remain unchanged. This preserves the stochastic structure of the
original problem while providing the spatial regularity needed in the
approximation argument.

All convolution estimates used below have already been established in
Section~\ref{subsec:kernel}.

\subsection{Galerkin construction of the viscous solution}
\label{subsec:wellposed}

Let
\[
H:=L^2(\R),
\qquad
V:=H^1(\R),
\]
so that
\[
V\hookrightarrow H\hookrightarrow V^*
\]
is a Gelfand triple.

Since negative Sobolev norms will be used in the time-regularity
estimates, the Galerkin basis must be compatible with an appropriate
Sobolev scale. We therefore use the Hermite basis. Let
\[
A:=-\partial_{xx}+x^2+1
\]
be the shifted harmonic oscillator and let
$\{e_k\}_{k\ge1}\subset\mathcal S(\R)$ denote its normalized
eigenfunctions,
\[
Ae_k=2k\,e_k.
\]
For $r\in\R$, define
\begin{equation}\label{eq:hermite-scale}
\mathcal S_r
:=
\left\{
u:
\|u\|_r^2
:=
\sum_{k\ge1}
(2k)^r|(u,e_k)|^2<\infty
\right\}.
\end{equation}
Then
\[
\mathcal S_0=H,
\qquad
\mathcal S_{-r}=(\mathcal S_r)^*,
\]
and
\begin{equation}\label{eq:hermite-sobolev}
\mathcal S_1
=
H^1(\R)\cap L^2(x^2\,\d x)
\hookrightarrow H^1(\R),
\qquad
H^{-1}(\R)\hookrightarrow\mathcal S_{-1}.
\end{equation}

For
\[
H_m:=\operatorname{span}\{e_1,\ldots,e_m\},
\]
let
\[
P_m u
:=
\sum_{k=1}^m(u,e_k)e_k .
\]
Since $P_m$ is a spectral projection of $A$,
\begin{equation}\label{eq:Pm-stable}
\|P_mu\|_r
\le
\|u\|_r,
\qquad
P_mu\longrightarrow u
\quad\text{in }\mathcal S_r,
\qquad
r\in\R .
\end{equation}
This stability is the main reason for choosing the Hermite basis: it
allows $P_m$ to be applied safely also in negative regularity spaces.

\subsubsection*{Galerkin approximation}

We seek an approximate solution of the form
\begin{equation}\label{eq:GalerkinAnsatz}
u_m^\varepsilon(t,x)
=
\sum_{k=1}^m
g_k^m(t)e_k(x),
\end{equation}
with initial datum
\[
u_m^\varepsilon(0)
=
P_mu_0^\varepsilon.
\]
The Galerkin approximation satisfies
\begin{align}
(\d u_m^\varepsilon,\varphi)
+
\big(
\partial_x
(f(u_m^\varepsilon)\mathcal U[u_m^\varepsilon]),
\varphi
\big)\,\d t
&=
-\varepsilon
(\partial_xu_m^\varepsilon,\partial_x\varphi)\,\d t +
(\Phi(u_m^\varepsilon),\varphi)\,\d W(t)
\nonumber\\
&\quad
+
\int_Z
(G(u_m^\varepsilon,z),\varphi)
\,\widetilde\eta(\d z,\d t)
\label{eq:GalerkinWeak}
\end{align}
for every $\varphi\in H_m$.

Choosing $\varphi=e_i$, $i=1,\ldots,m$, gives
\begin{align}
\d g_i^m
&=
-
\big(
\partial_x
(f(u_m^\varepsilon)\mathcal U[u_m^\varepsilon]),
e_i
\big)\,\d t
-
\varepsilon
(\partial_xu_m^\varepsilon,\partial_xe_i)\,\d t
\nonumber\\
&\quad
+
(\Phi(u_m^\varepsilon),e_i)\,\d W(t)
+
\int_Z
(G(u_m^\varepsilon,z),e_i)
\,\widetilde\eta(\d z,\d t).
\label{eq:GalerkinCoeff}
\end{align}

Equivalently,
\begin{equation}\label{eq:VectorSDE}
\mathrm{d}\mathbf g^m
=
F_m(\mathbf g^m)\,\d t
+
B_m(\mathbf g^m)\,\d W
+
\int_Z
J_m(\mathbf g^m_-,z)
\,\widetilde\eta(\d z,\d t).
\end{equation}

\begin{lemma}\label{lem:coefficients}
For every fixed $m$, the mappings $F_m$ and $B_m$ are locally
Lipschitz on $\R^m$, and
\[
\int_Z
|J_m(\mathbf a,z)-J_m(\mathbf b,z)|^2
\,\nu(\d z)
\le
C_m|\mathbf a-\mathbf b|^2.
\]
Moreover,
\[
|F_m(\mathbf a)|
\le
C_m(1+|\mathbf a|^2),
\]
whereas
\[
|B_m(\mathbf a)|
+
\left(
\int_Z
|J_m(\mathbf a,z)|^2
\,\nu(\d z)
\right)^{1/2}
\le
C_m|\mathbf a|.
\]
\end{lemma}

\begin{proof}
The local Lipschitz estimates follow from the global Lipschitz
continuity of $f,\Phi$ and $G$, Young's convolution inequality, and
equivalence of norms on the finite-dimensional space $H_m$. For the
transport term, we use
\[
f(u_{\mathbf a})\mathcal U[u_{\mathbf a}]
-
f(u_{\mathbf b})\mathcal U[u_{\mathbf b}]
=
(f(u_{\mathbf a})-f(u_{\mathbf b}))
\mathcal U[u_{\mathbf a}]
+
f(u_{\mathbf b})
\big(
\mathcal U[u_{\mathbf a}]
-
\mathcal U[u_{\mathbf b}]
\big).
\]
Since both $f(u_{\mathbf a})$ and
$\mathcal U[u_{\mathbf a}]$ have at most linear growth in
$\mathbf a$, their product has at most quadratic growth. Thus no
linear-growth estimate is claimed for $F_m$.
\end{proof}

The standard existence theory for finite-dimensional stochastic
differential equations with jumps therefore yields a unique maximal
adapted càdlàg solution, defined up to an explosion time $\tau_m$.

\subsubsection{Stopped energy estimate}
\label{subsubsec:truncation}

For $R>0$, define
\[
\tau_R^m
:=
\inf
\left\{
t\ge0:
\|u_m^\varepsilon(t)\|_{L^2}\ge R
\right\}
\wedge T.
\]
Before $\tau_R^m$,
\[
\|u_m^\varepsilon(t)\|_{L^2}\le R,
\]
and therefore, by Young's inequality,
\begin{equation}\label{eq:stopped-U-bounds}
\|\mathcal U[u_m^\varepsilon]\|_{L^\infty}
\le
R\|\mu\|_{L^2},
\qquad
\|\partial_x\mathcal U[u_m^\varepsilon]\|_{L^\infty}
\le
R\|\mu'\|_{L^2}.
\end{equation}
Thus the stopping time provides precisely the deterministic control on
the nonlocal coefficient required for the energy estimate.

\begin{remark}[Why the $L^1$-estimate is postponed]
\label{rem:circularity}
The natural $L^1$ estimate uses the nonlinear test function
$\beta_\delta'(u)$, while the $BV$ estimate uses
$\beta_\delta'(\partial_xu)$. These functions do not belong to $H_m$
in general. Testing the Galerkin equation would instead introduce
their projections $P_m\beta_\delta'(u_m^\varepsilon)$, which cannot be
replaced by the original nonlinear functions because $P_m$ does not
commute with nonlinear superposition. We therefore derive these
estimates only after passing to the viscous equation.
\end{remark}

\begin{proposition}[Stopped energy estimate]
\label{prop:energy-global-galerkin}
For every fixed $\varepsilon>0$ and $R>0$, there exists
$C_R>0$, independent of $m$, such that
\begin{align}
\mathbb E
\sup_{0\le t\le T}
\|u_m^\varepsilon(t\wedge\tau_R^m)\|_{L^2}^2
+
\varepsilon
\mathbb E
\int_0^{T\wedge\tau_R^m}
\|\partial_xu_m^\varepsilon(t)\|_{L^2}^2\,\d t
\le
C_R .
\label{eq:global-energy-final}
\end{align}
\end{proposition}

\begin{proof}
Applying the finite-dimensional It\^o--L\'evy formula to
$\|u_m^\varepsilon\|_{L^2}^2$ up to $\tau_R^m$ gives
\begin{align*}
\|u_m^\varepsilon(t\wedge\tau_R^m)\|_{L^2}^2
&+
2\varepsilon
\int_0^{t\wedge\tau_R^m}
\|\partial_xu_m^\varepsilon\|_{L^2}^2\,\d s =
\|P_mu_0^\varepsilon\|_{L^2}^2
+
\mathcal I_{\rm tr}(t)
+
\mathcal I_{\rm corr}(t)
+
M_B^R(t)+M_J^R(t),
\end{align*}
where $\mathcal I_{\rm tr}$ denotes the transport contribution and
$\mathcal I_{\rm corr}$ contains the Brownian and jump It\^o
corrections.

Set
\[
F_0(r):=\int_0^r f(s)\,\d s.
\]
Since $f(0)=0$ and $f$ is $L_f$-Lipschitz, $|F_0(r)|
\le
\frac{L_f}{2}|r|^2.$ Integration by parts gives
\[
-\big(
u_m^\varepsilon,
\partial_x(f(u_m^\varepsilon)
\mathcal U[u_m^\varepsilon])
\big)
=
-\int_\R
F_0(u_m^\varepsilon)
\partial_x\mathcal U[u_m^\varepsilon]\,\d x.
\]
Hence, by \eqref{eq:stopped-U-bounds},
\[
|\mathcal I_{\rm tr}'(t)|
\le
C_R
\|u_m^\varepsilon(t)\|_{L^2}^2 .
\]

The homogeneous noise assumptions imply
\[
\|\Phi(u)\|_{L^2}^2
+
\int_Z
\|G(u,z)\|_{L^2}^2\,\nu(\d z)
\le
C\|u\|_{L^2}^2.
\]
For the Brownian martingale, BDG and Young's inequalities yield
\[
\mathbb E
\sup_{s\le t}|M_B^R(s)|
\le
\theta
\mathbb E
\sup_{s\le t}
\|u_m^\varepsilon(s\wedge\tau_R^m)\|_{L^2}^2
+
C_\theta
\int_0^t
\mathbb E
\|u_m^\varepsilon(s\wedge\tau_R^m)\|_{L^2}^2\,\d s .
\]

For the compensated Poisson martingale, write
\[
H(u,z)
:=
\|u+G(u,z)\|_{L^2}^2-\|u\|_{L^2}^2
=
2(u,G(u,z))_{L^2}
+
\|G(u,z)\|_{L^2}^2.
\]
Using BDG with exponent one and
\[
\int_Z|H(u,z)|\,\nu(\d z)
\le
C\|u\|_{L^2}^2,
\]
we obtain
\[
\mathbb E
\sup_{s\le t}|M_J^R(s)|
\le
C
\mathbb E
\int_0^{t\wedge\tau_R^m}
\|u_m^\varepsilon(s)\|_{L^2}^2\,\d s.
\]
This argument uses only the $L^1(Z,\nu)$- and $L^2(Z,\nu)$-bounds of
\textbf{(A3)} and therefore requires no fourth moment of $G$.

Combining the estimates, absorbing the Brownian supremum term, and
applying Gr\"onwall's inequality proves
\eqref{eq:global-energy-final}.
\end{proof}

\subsubsection{Stopped pathwise uniqueness}

\begin{lemma}[Stopped pathwise uniqueness]
\label{lem:pathwise-uniq}
Let $u$ and $v$ be two solutions of the viscous equation defined on the
same stochastic basis, driven by the same $W$ and $\eta$, and having
the same initial datum. For $R>0$, define
\[
\rho_R
:=
\inf
\left\{
t\ge0:
\|u(t)\|_{L^2}
+
\|v(t)\|_{L^2}
\ge R
\right\}
\wedge T .
\]
Then
\[
u(t)=v(t),
\qquad
0\le t\le\rho_R,
\quad\text{a.s.}
\]
\end{lemma}

\begin{proof}
Set
\[
w:=u-v \qquad \mbox{and} \qquad D :=
f(u)\mathcal U[u] - f(v)\mathcal U[v]\,.
\]

On $[0,\rho_R]$,
\begin{align*}
D
&=
(f(u)-f(v))\mathcal U[u]
+
f(v)\big(\mathcal U[u]-\mathcal U[v]\big),
\end{align*}
and the convolution estimates imply
\[
\|D\|_{L^2}
\le
C_R\|w\|_{L^2}.
\]
Applying It\^o's formula to $\|w\|_{L^2}^2$, using
\textbf{(A2)}--\textbf{(A3)} and Young's inequality, gives
\[
\d\|w\|_{L^2}^2
+
\varepsilon
\|\partial_xw\|_{L^2}^2\,\d t
\le
C_{R,\varepsilon}
\|w\|_{L^2}^2\,\d t
+
\d N_t.
\]
Stopping, taking expectations, and applying Gr\"onwall's inequality
give
\[
\mathbb E
\|w(t\wedge\rho_R)\|_{L^2}^2
=0,
\]
since $w(0)=0$.
\end{proof}

\subsection{Compactness of the stopped Galerkin approximations}
\label{subsec:galerkin-compactness}
\label{subsec:compactness}

For fixed $\varepsilon,R>0$, set
\[
u_{m,R}^\varepsilon(t)
:=
u_m^\varepsilon(t\wedge\tau_R^m).
\]
Besides the spatial energy estimate, we need temporal regularity in a
negative Sobolev space.

\begin{proposition}[Time-translation estimate]
\label{prop:time-translation}
For every fixed $\varepsilon,R>0$,
\begin{equation}\label{eq:time-transl}
\sup_{m\ge1}
\mathbb E
\int_0^{T-h}
\|
u_{m,R}^\varepsilon(t+h)
-
u_{m,R}^\varepsilon(t)
\|_{\mathcal S_{-1}}^2
\,\d t
\le
C_R(h+\varepsilon h),
\end{equation}
for every $0<h<T$. Consequently, for every
$\alpha\in(0,\tfrac12)$,
\begin{equation}\label{eq:frac-time}
\sup_{m\ge1}
\mathbb E
\|u_{m,R}^\varepsilon\|_{W^{\alpha,2}
(0,T;\mathcal S_{-1})}^2
<\infty .
\end{equation}
\end{proposition}

\begin{proof}
Integrating the Galerkin equation from $t$ to $t+h$ gives four
contributions corresponding to transport, viscosity, Brownian noise,
and jump noise.

By \eqref{eq:Pm-stable},
\[
\|P_mg\|_{\mathcal S_{-1}}
\le
\|g\|_{\mathcal S_{-1}}
\le
C\|g\|_{H^{-1}},
\]
so it suffices to estimate the increments in $H^{-1}(\R)$.

For the transport term,
\[
\|\partial_x q\|_{H^{-1}}
\le
\|q\|_{L^2},
\]
and on $[0,\tau_R^m]$,
\[
\|f(u_m^\varepsilon)
\mathcal U[u_m^\varepsilon]\|_{L^2}
\le
C_R\|u_m^\varepsilon\|_{L^2}.
\]
Using Cauchy--Schwarz in time and
\eqref{eq:global-energy-final}, this contribution is bounded by
$C_Rh$ after integration in $t$.

Similarly,
\[
\|\partial_{xx}u\|_{H^{-1}}
\le
\|\partial_xu\|_{L^2},
\]
so the viscous contribution is bounded by $C_R\varepsilon h$.
The It\^o and compensated-Poisson isometries, together with the
homogeneous growth assumptions on $\Phi$ and $G$, give an $O(h)$
bound for the stochastic increments.

Combining the four estimates proves \eqref{eq:time-transl}. Finally,
\[
\int_0^T
h^{-1-2\alpha}
\left(
\int_0^{T-h}
\|u(t+h)-u(t)\|_{\mathcal S_{-1}}^2\,\d t
\right)\d h
\]
is finite for every $\alpha<1/2$, which yields
\eqref{eq:frac-time}.
\end{proof}

\begin{remark}
The fractional estimate \eqref{eq:frac-time}, rather than a supremum
over all time increments, is the compactness quantity used below. In
particular, one cannot interchange expectation with
\[
\sup_{0<h<T}
h^{-1}
\int_0^{T-h}
\|u(t+h)-u(t)\|_{\mathcal S_{-1}}^2\,\d t .
\]
\end{remark}

\begin{proposition}[Compactness of the stopped Galerkin sequence]
\label{prop:compactness}
For every fixed $\varepsilon,R>0$, there exists a progressively
measurable process $u_R^\varepsilon$ such that, for every bounded
interval $I\subset\R$,
\begin{align}
u_{m,R}^\varepsilon
&\rightharpoonup
u_R^\varepsilon
\quad\text{weakly in }
L^2\big(
\Omega;L^2(0,T;H^1(I))
\big),
\label{eq:weak-H1}
\\
u_{m,R}^\varepsilon
&\longrightarrow
u_R^\varepsilon
\quad\text{in probability in }
L^2(0,T;L^2(I)),
\label{eq:strong-L2}
\\
u_{m,R}^\varepsilon
&\longrightarrow
u_R^\varepsilon
\quad\text{strongly in }
L^q\big(
\Omega;L^2(0,T;L^2(I))
\big),
\qquad q<2.
\label{eq:strong-Lq}
\end{align}
Moreover, along a subsequence,
\begin{equation}\label{eq:ae}
u_{m,R}^\varepsilon(t,x)
\longrightarrow
u_R^\varepsilon(t,x)
\qquad
\text{for a.e. }(\omega,t,x),
\end{equation}
and
\[
u_R^\varepsilon
\in
D([0,T];\mathcal S_{-1})
\qquad\mathbb P\text{-a.s.}
\]
\end{proposition}

\begin{proof}
Let $I\subset\R$ be bounded. Since
\[
H^1(I)
\Subset
L^2(I)
\hookrightarrow
\mathcal S_{-1},
\]
the compact embedding
\begin{equation}\label{eq:frac-compact}
L^2(0,T;H^1(I))
\cap
W^{\alpha,2}(0,T;\mathcal S_{-1})
\Subset
L^2(0,T;L^2(I))
\end{equation}
holds for $\alpha\in(0,\tfrac12)$. Consequently,
\eqref{eq:global-energy-final}, \eqref{eq:frac-time}, and Markov's
inequality imply tightness of the laws of
$\{u_{m,R}^\varepsilon\}$ in the local space-time topology.

The increment estimate also yields the Aldous condition in
$\mathcal S_{-1}$, and therefore tightness in
$D([0,T];\mathcal S_{-1})$.

For the noise variables we use their natural canonical spaces:
\[
W\in C([0,T];\R),
\]
while the Poisson noise is represented by the integer-valued random
measure
\[
\eta\in\mathcal M_{\mathbb N}([0,T]\times Z).
\]
The compensated measure
\[
\widetilde\eta
=
\eta-\nu(\d z)\d t
\]
is reconstructed after passage to the limit. Thus, $\left\{
\mathcal L
\big(
u_{m,R}^\varepsilon,W,\eta
\big)
\right\}_{m\ge1}$ is tight.

After applying a Skorokhod representation, let
$\widetilde u_m\to\widetilde u$ almost surely in the resulting
topology. The limit equation is identified as follows.

\smallskip
\noindent
\emph{Admissible test functions.}
The Galerkin equation is valid only for test functions in $H_m$.
Given $\varphi\in\mathcal S(\R)$, we therefore use $\varphi_m:=P_m\varphi.$ By \eqref{eq:Pm-stable}, $\varphi_m\longrightarrow\varphi$ in every $\mathcal S_r$, and hence in all the norms needed to pass to
the deterministic terms.

\smallskip
\noindent
\emph{The nonlocal term.}
Since only local strong convergence is available, the global Young
estimate cannot be applied directly to
$\widetilde u_m-\widetilde u$. We write
\[
\mu=\mu_L+\mu^L,
\qquad
\mu_L:=\mu\mathbf1_{[-L,L]}.
\]
On compact spatial sets, the convolution with $\mu_L$ depends only on
the solution on a larger compact set, so local strong convergence
applies. For the tail,
\[
\|\mu^L*(\widetilde u_m-\widetilde u)\|_{L^2}
\le
\|\mu^L\|_{L^1}
\|\widetilde u_m-\widetilde u\|_{L^2},
\]
and the right-hand side is uniformly small as $L\uparrow\infty$ by
the global $L^2$-bounds.

\smallskip
\noindent
\emph{The stochastic integrals.}
They are not identified by continuity of the It\^o map with respect
to the Skorokhod convergence. Instead, one passes to the limit in the
martingale problem. For $\varphi\in\mathcal S(\R)$, define
\[
M_\varphi(t)
:=
(\widetilde u(t),\varphi)
-
(u_0^\varepsilon,\varphi)
-
\int_0^t
\Big[
(f(\widetilde u)\mathcal U[\widetilde u],
\partial_x\varphi)
-
\varepsilon
(\partial_x\widetilde u,\partial_x\varphi)
\Big]\d s .
\]
Passing to the limit in the Galerkin martingales identifies
$M_\varphi$ as a square-integrable martingale with Brownian quadratic
variation determined by $\Phi(\widetilde u)$ and jump compensator
determined by $G(\widetilde u,z)$. The martingale representation
theorem then recovers the required Brownian and Poisson stochastic
integrals.

Thus every subsequential limit solves the same stopped viscous problem.
By Lemma~\ref{lem:pathwise-uniq}, two such limits driven by the same
noise coincide. The Gy\"ongy--Krylov criterion therefore implies
convergence in probability of the full Galerkin sequence on the
original probability space. Finally, the uniform $L^2$-moment bounds
and Vitali's theorem imply \eqref{eq:strong-Lq}. The almost-everywhere
convergence follows after extracting a subsequence.
\end{proof}

\subsubsection*{Identification of the Galerkin limit}

\begin{proposition}[Identification of the Galerkin limit]
\label{prop:Galerkin-limit}
For every $\varphi\in C_c^\infty(\R)$,
$u_R^\varepsilon$ satisfies
\begin{align}
(u_R^\varepsilon(t),\varphi)
&=
(u_0^\varepsilon,\varphi)
+
\int_0^{t\wedge\tau_R^\varepsilon}
\big(
f(u_R^\varepsilon)
\mathcal U[u_R^\varepsilon],
\partial_x\varphi
\big)\,\d s
\nonumber\\
&\quad
-
\varepsilon
\int_0^{t\wedge\tau_R^\varepsilon}
(\partial_xu_R^\varepsilon,\partial_x\varphi)\,\d s
\nonumber\\
&\quad
+
\int_0^{t\wedge\tau_R^\varepsilon}
(\Phi(u_R^\varepsilon),\varphi)\,\d W(s)
\nonumber\\
&\quad
+
\int_0^{t\wedge\tau_R^\varepsilon}
\!\!\int_Z
(G(u_R^\varepsilon,z),\varphi)
\,\widetilde\eta(\d z,\d s).
\label{eq:viscous-limit}
\end{align}
Hence $u_R^\varepsilon$ solves the viscous equation up to its
$L^2$-exit time.
\end{proposition}

\begin{proof}
The deterministic and stochastic passages to the limit have already
been identified in the proof of Proposition~\ref{prop:compactness}.
The present proposition merely records the resulting weak formulation.
\end{proof}

By stopped pathwise uniqueness, the family
$\{u_R^\varepsilon\}_{R\ge1}$ is consistent. Hence these solutions can
be patched together to obtain a maximal viscous solution
$u^\varepsilon$ on
\[
[0,\tau_\infty^\varepsilon),
\qquad
\tau_R^\varepsilon
\uparrow
\tau_\infty^\varepsilon.
\]

\subsection{A priori estimates for the viscous solution}
\label{subsec:viscous-estimates}

We now work directly with the viscous solution. At this level the
equation is no longer projected, so the entropy and variation test
functions can be used legitimately.

\subsubsection{Uniform $L^1$-estimate}

\begin{proposition}[Uniform $L^1$-estimate]
\label{prop:L1}
The maximal viscous solution satisfies
\begin{align}
\mathbb E
\sup_{t<\tau_\infty^\varepsilon\wedge T}
\|u^\varepsilon(t)\|_{L^1}
&\le
C
\left(
1+\mathbb E\|u_0\|_{L^1}
\right),
\label{eq:L1-uniform}
\\
\mathbb E
\sup_{t<\tau_\infty^\varepsilon\wedge T}
\|u^\varepsilon(t)\|_{L^1}^2
&\le
C
\left(
1+\mathbb E\|u_0\|_{L^1}^2
\right),
\label{eq:L1-second-moment}
\end{align}
where $C$ is independent of $\varepsilon$.
\end{proposition}

\begin{proof}
Fix an $L^2$-exit time $\tau_R^\varepsilon$ and introduce a spatial
cut-off $\chi_L$ satisfying
\[
0\le\chi_L\le1,
\qquad
\chi_L=1\ \text{on }[-L,L],
\qquad
|\chi_L'|\le C/L.
\]
Let
\[
X_{\delta,L}(t)
:=
\int_\R
\beta_\delta(u^\varepsilon(t,x))
\chi_L(x)\,\d x.
\]

Applying the It\^o--L\'evy formula and setting
\[
\Xi(r)
:=
\int_0^r
\beta_\delta'(s)f'(s)\,\d s,
\]
the transport contribution can be written as
\begin{align*}
\mathcal T_{\delta,L}
&=
\int_\R
\Big[
\Xi(u^\varepsilon)
-
\beta_\delta'(u^\varepsilon)f(u^\varepsilon)
\Big]
\partial_x\mathcal U[u^\varepsilon]
\chi_L\,\d x +
\int_\R
\Xi(u^\varepsilon)
\mathcal U[u^\varepsilon]
\chi_L'\,\d x .
\end{align*}

The identity
\begin{equation}\label{eq:L1-transport-identity}
\Xi(r)
-
\beta_\delta'(r)f(r)
=
-\int_0^r
\beta_\delta''(s)f(s)\,\d s
\end{equation}
and $|f(s)|\le L_f|s|$ imply
\[
\left|
\Xi(r)
-
\beta_\delta'(r)f(r)
\right|
\le
2L_f\delta .
\]
Thus the first transport term vanishes as $\delta\downarrow0$ for
fixed $R$ and $L$.

For the boundary term, $|\Xi(r)|\le L_f|r|,$ and on $[0,\tau_R^\varepsilon]$,
\[
\|\mathcal U[u^\varepsilon]\|_\infty
\le
R\|\mu\|_2.
\]
Hence
\[
\left|
\int_\R
\Xi(u^\varepsilon)
\mathcal U[u^\varepsilon]\chi_L'\,\d x
\right|
\le
C_R
\|u^\varepsilon\|_2
\|\chi_L'\|_2
\le
C_RL^{-1/2},
\]
which vanishes as $L\uparrow\infty$.

The viscosity term is nonpositive. Moreover,
\[
\frac12
\int_\R
\beta_\delta''(u^\varepsilon)
\Phi(u^\varepsilon)^2\,\d x
\le
\frac12L_\Phi^2
\|u^\varepsilon\|_{L^1},
\]
because $\beta_\delta''(r)r^2\le|r|.$ Likewise, by Lemma~\ref{lem:jump-convex},
\[
0\le
\int_\R\!\!\int_Z
\Big[
\beta_\delta(u^\varepsilon+G)
-
\beta_\delta(u^\varepsilon)
-
\beta_\delta'(u^\varepsilon)G
\Big]
\nu(\d z)\,\d x
\le
2\widetilde L_G
\|u^\varepsilon\|_{L^1}.
\]

The Brownian and jump martingales are estimated by BDG. Passing first
$\delta\downarrow0$ and then $L\uparrow\infty$, followed by
Gr\"onwall's inequality, yields
\eqref{eq:L1-uniform}. The resulting estimate is independent of the
$L^2$-exit level $R$, so $R\uparrow\infty$ gives the estimate up to
$\tau_\infty^\varepsilon$.

For the second moment, square the corresponding pathwise estimate.
The Brownian integrand satisfies
\[
\left|
\int_\R
\beta_\delta'(u)\Phi(u)\,\d x
\right|
\le
L_\Phi\|u\|_{L^1},
\]
while Minkowski's integral inequality and \textbf{(A3)} give, for the
Poisson integrand,
\[
\left(
\int_Z
\left|
\int_\R
[
\beta_\delta(u+G)-\beta_\delta(u)
]\,\d x
\right|^2
\nu(\d z)
\right)^{1/2}
\le
L_G\|u\|_{L^1}.
\]
A second application of BDG and Gr\"onwall gives
\eqref{eq:L1-second-moment}.
\end{proof}

\subsubsection{Global continuation of the viscous solution}

\begin{proposition}[Global viscous solution]
\label{prop:remove-truncation}
For every $\varepsilon>0$,
\[
\tau_\infty^\varepsilon=T
\qquad\text{a.s.}
\]
Consequently, \eqref{eq:viscous} admits a unique global viscous
solution satisfying
\begin{align}
\mathbb E
\sup_{0\le t\le T}
\|u^\varepsilon(t)\|_{L^1}
&\le
C
\left(
1+\mathbb E\|u_0\|_{L^1}
\right),
\label{eq:L1-untruncated}
\\
\mathbb E
\sup_{0\le t\le T}
\|u^\varepsilon(t)\|_{L^1}^2
&\le
C
\left(
1+\mathbb E\|u_0\|_{L^1}^2
\right).
\label{eq:L1-untruncated-2}
\end{align}
\end{proposition}

\begin{proof}
For $N\ge1$, define
\[
\sigma_N^\varepsilon
:=
\inf
\left\{
t\ge0:
\|u^\varepsilon(t)\|_{L^1}\ge N
\right\}
\wedge T.
\]
On $[0,\sigma_N^\varepsilon]$,
\[
\|\partial_x\mathcal U[u^\varepsilon]\|_{L^\infty}
\le
N\|\mu'\|_{L^\infty},
\]
which is now a deterministic bound. Repeating the energy argument up
to $\sigma_N^\varepsilon\wedge\tau_R^\varepsilon$ therefore yields
\begin{equation}\label{eq:energy-localised}
\mathbb E
\sup_{t\le
\sigma_N^\varepsilon\wedge\tau_R^\varepsilon}
\|u^\varepsilon(t)\|_{L^2}^2
\le
C_N ,
\end{equation}
where $C_N$ is independent of $R$ and $\varepsilon$.

By \eqref{eq:L1-uniform} and Markov's inequality,
\[
\mathbb P(\sigma_N^\varepsilon<T)
\le
\frac{C}{N}.
\]
Consequently,
\[
\mathbb P(\tau_R^\varepsilon<T)
\le
\mathbb P(\sigma_N^\varepsilon<T)
+
\mathbb P\left(
\sup_{t\le
\sigma_N^\varepsilon\wedge\tau_R^\varepsilon}
\|u^\varepsilon(t)\|_{L^2}
\ge R
\right)
\le
\frac{C}{N}
+
\frac{C_N}{R^2}.
\]
Given $\eta>0$, choose first $N$ and then $R$ sufficiently large.
Thus
\[
\mathbb P(\tau_\infty^\varepsilon<T)=0.
\]
The estimates
\eqref{eq:L1-untruncated}--\eqref{eq:L1-untruncated-2}
then follow from Proposition~\ref{prop:L1}.
\end{proof}

\subsubsection{Localized $BV$-estimate}
\label{subsubsec:BV-estimate}

We now impose \textbf{(K4)} in addition to the standing assumptions.
The proof below uses spatial translations, and hence only the Lipschitz
regularity in \textbf{(A1)}--\textbf{(A3)}; in particular, no second
derivative of $f$ is required.

\begin{proposition}[$BV$-estimate]
\label{prop:BV}
Let
\begin{equation}\label{eq:BV-stopping}
\tau_M^\varepsilon
:=
\inf
\left\{
t\ge0:
\|u^\varepsilon(t)\|_{L^1}
+
\|u^\varepsilon(t)\|_{L^2}
\ge M
\right\}
\wedge T.
\end{equation}
Then, for every $M\ge1$,
\begin{equation}\label{eq:BV-stopped}
\mathbb E
\sup_{0\le t\le\tau_M^\varepsilon}
\TV(u^\varepsilon(t))
\le
C_M
\left(
1+\mathbb E\TV(u_0)
\right),
\end{equation}
where $C_M$ is independent of $\varepsilon$. Furthermore,
\begin{equation}\label{eq:tau-prob}
\lim_{M\uparrow\infty}
\sup_{\varepsilon>0}
\mathbb P(
\tau_M^\varepsilon<T
)
=
0,
\end{equation}
and consequently
\begin{equation}\label{eq:TV-bdd-in-prob}
\lim_{\Lambda\uparrow\infty}
\sup_{\varepsilon>0}
\mathbb P
\left(
\sup_{0\le t\le T}
\TV(u^\varepsilon(t))
\ge\Lambda
\right)
=
0.
\end{equation}
\end{proposition}

\begin{proof}
For $h\neq0$, denote the spatial translation by
$\tau_h w(x):=w(x+h)$ and set
\[
V_h(w)
:=
\frac{1}{|h|}
\|\tau_h w-w\|_{L^1(\R)}.
\]
Recall that, for every $w\in BV(\R)$,
\begin{equation}\label{eq:BV-dq-char}
V_h(w)\le \TV(w),
\qquad
\lim_{h\to0}V_h(w)=\TV(w).
\end{equation}

Fix $\varepsilon>0$ and write
\[
u:=u^\varepsilon,
\qquad
u_h:=\tau_hu,
\qquad
d_h:=u_h-u,
\qquad
\mathcal U:=\mathcal U[u],
\qquad
\mathcal U_h:=\tau_h\mathcal U.
\]
By translation invariance of the convolution,
$\mathcal U_h=\mathcal U[u_h]$. Subtracting the viscous equation for
$u$ from its spatial translate gives
\begin{align}\label{eq:translated-difference}
\d d_h
&+
\partial_x\!\left(
 f(u_h)\mathcal U_h-f(u)\mathcal U
\right)\d t
=
\varepsilon\partial_{xx}d_h\,\d t
+
\big(\Phi(u_h)-\Phi(u)\big)\,\d W(t)
\nonumber\\
&\quad
+
\int_Z
\big(G(u_h,z)-G(u,z)\big)
\widetilde\eta(\d z,\d t).
\end{align}
We apply the It\^o--L\'evy formula to
$\int_\R\beta_\delta(d_h)\,\d x$, first with a compact spatial cut-off
and then let the cut-off tend to one. All estimates below are independent
of the cut-off.

For the transport term, decompose
\begin{equation}\label{eq:BV-flux-split}
f(u_h)\mathcal U_h-f(u)\mathcal U
=
\big(f(u_h)-f(u)\big)\mathcal U_h
+
f(u)\big(\mathcal U_h-\mathcal U\big).
\end{equation}
The first term in \eqref{eq:BV-flux-split} is the usual local
conservation-law contribution. Indeed, with
\[
q_h
:=
\sgn(d_h)\big(f(u_h)-f(u)\big),
\]
the Kato chain rule for a Lipschitz flux gives, in the sense of
distributions,
\[
\sgn(d_h)\,
\partial_x\big(f(u_h)-f(u)\big)
=
\partial_x q_h.
\]
Consequently, after letting $\delta\downarrow0$,
\[
-\int_\R
\sgn(d_h)
\partial_x\!\left(
\big(f(u_h)-f(u)\big)\mathcal U_h
\right)\d x
=
-\int_\R
\partial_x\big(\mathcal U_hq_h\big)\,\d x
=0.
\]
Thus no derivative of $f'$ is involved.

It remains to estimate the genuinely nonlocal remainder
$f(u)(\mathcal U_h-\mathcal U)$. Since $f$ is Lipschitz and $f(0)=0$,
the Sobolev chain rule gives
$|\partial_x f(u)|\le L_f|\partial_xu|$ a.e. Hence
\begin{align}
&\left|
\int_\R
\sgn(d_h)
\partial_x\!\left(
 f(u)(\mathcal U_h-\mathcal U)
\right)\d x
\right|\le
L_f
\|\mathcal U_h-\mathcal U\|_{L^\infty}
\TV(u)
+
L_f
\int_\R
|u|\,
|\partial_x\mathcal U_h-\partial_x\mathcal U|\,\d x.
\label{eq:BV-nonlocal-dq}
\end{align}
By \textbf{(K3)}--\textbf{(K4)},
\[
\|\mathcal U_h-\mathcal U\|_{L^\infty}
\le
|h|\,\|\partial_x\mathcal U\|_{L^\infty}
\le
|h|\,\|\mu'\|_{L^\infty}\|u\|_{L^1},
\]
and, since
$\partial_x\mathcal U=\mu'*u$,
\begin{align*}
\|\partial_x\mathcal U_h-\partial_x\mathcal U\|_{L^2}
&\le
\|\tau_h\mu'-\mu'\|_{L^1}\|u\|_{L^2}\le
|h|\,\|\mu''\|_{L^1}\|u\|_{L^2}.
\end{align*}
Therefore, after division by $|h|$ and passage to $h\to0$,
the nonlocal transport contribution is bounded by
\begin{equation}\label{eq:BV-T2-bound}
L_f\|\mu'\|_{L^\infty}
\|u\|_{L^1}\TV(u)
+
L_f\|\mu''\|_{L^1}\|u\|_{L^2}^2.
\end{equation}
In particular, on $[0,\tau_M^\varepsilon]$ it is bounded by
$C_M(1+\TV(u))$.

The viscosity contribution is nonpositive. For the Brownian It\^o
correction, \textbf{(A2)} and
$\beta_\delta''(r)r^2\le |r|$ give
\[
\frac12
\int_\R
\beta_\delta''(d_h)
|\Phi(u_h)-\Phi(u)|^2\,\d x
\le
\frac12L_\Phi^2\|d_h\|_{L^1}.
\]
For the jump compensator, using
\[
0\le
\beta_\delta(a+b)-\beta_\delta(a)-\beta_\delta'(a)b
\le 2|b|
\]
and the $L^1(Z,\nu)$-Lipschitz bound in \textbf{(A3)}, we obtain
\[
\int_Z\!\int_\R
\Big[
\beta_\delta(d_h+\Delta_hG)
-
\beta_\delta(d_h)
-
\beta_\delta'(d_h)\Delta_hG
\Big] \,\d x\,\nu(\d z)
\le
2\widetilde L_G\|d_h\|_{L^1},
\]
where
$\Delta_hG:=G(u_h,z)-G(u,z)$.

The martingale terms are controlled directly by the Lipschitz assumptions.
After division by $|h|$, the Brownian integrand satisfies
\[
\frac1{|h|}
\left|
\int_\R
\beta_\delta'(d_h)
\big(\Phi(u_h)-\Phi(u)\big)\,\d x
\right|
\le
L_\Phi V_h(u),
\]
while Minkowski's inequality and the $L^2(Z,\nu)$-Lipschitz bound in
\textbf{(A3)} give the analogous estimate
\[
\frac1{|h|}
\left(
\int_Z
\left|
\int_\R
\big[
\beta_\delta(d_h+\Delta_hG)-\beta_\delta(d_h)
\big]\,\d x
\right|^2
\nu(\d z)
\right)^{1/2}
\le
L_G V_h(u).
\]
Hence the BDG inequality yields, for every $\theta\in(0,1)$,
\begin{align}\label{eq:BV-dq-BDG}
\mathbb E\sup_{s\le t\wedge\tau_M^\varepsilon}|N_h(s)|
&\le
\theta\,
\mathbb E\sup_{s\le t\wedge\tau_M^\varepsilon}V_h(u(s))+
C_\theta
\int_0^t
\mathbb E\sup_{r\le s\wedge\tau_M^\varepsilon}V_h(u(r))\,\d s,
\end{align}
where $N_h$ denotes the sum of the normalized Brownian and Poisson
martingales.

For complete rigor in the limit $h\to0$, one may first stop additionally at
\[
\rho_K^\varepsilon
:=
\inf\{t\ge0:\TV(u^\varepsilon(t))\ge K\}\wedge T.
\]
On $[0,\tau_M^\varepsilon\wedge\rho_K^\varepsilon]$ all variation terms
are integrable. We first let $\delta\downarrow0$ and then $h\to0$ in the
preceding estimates, using \eqref{eq:BV-dq-char}, Fatou's lemma, and the
lower semicontinuity of the total variation. Since, by Lemma~\ref{lem:init-approx},
\[
V_h(u_0^\varepsilon)
\le
\TV(u_0^\varepsilon)
\le
\TV(u_0)
\qquad\text{a.s.},
\]
we obtain, after absorbing the $\theta$-term,
\begin{align*}
\mathbb E
\sup_{s\le t\wedge\tau_M^\varepsilon\wedge\rho_K^\varepsilon}
\TV(u^\varepsilon(s))
&\le
\mathbb E\TV(u_0)
+
C_M
\int_0^t
\left[
1+
\mathbb E
\sup_{r\le s\wedge\tau_M^\varepsilon\wedge\rho_K^\varepsilon}
\TV(u^\varepsilon(r))
\right]\d s.
\end{align*}
Gr\"onwall's inequality gives a bound independent of $K$ and
$\varepsilon$. Letting $K\uparrow\infty$ proves
\eqref{eq:BV-stopped}.

It remains to estimate the exit probability. Define
\[
\varsigma_N^\varepsilon
:=
\inf
\left\{
t\ge0:
\|u^\varepsilon(t)\|_{L^1}\ge N
\right\}
\wedge T.
\]
On $[0,\varsigma_N^\varepsilon]$, the localized energy estimate gives
\begin{equation}\label{eq:energy-localised-eps}
\mathbb E
\sup_{t\le\varsigma_N^\varepsilon}
\|u^\varepsilon(t)\|_{L^2}^2
\le
C_N,
\end{equation}
uniformly in $\varepsilon$. Therefore, for $M>2N$,
\[
\mathbb P(\tau_M^\varepsilon<T)
\le
\frac{C}{N}
+
\frac{4C_N}{M^2}.
\]
Choosing first $N$ and then $M$ proves \eqref{eq:tau-prob}.

Finally,
\[
\mathbb P
\left(
\sup_{t\le T}\TV(u^\varepsilon(t))
\ge\Lambda
\right)
\le
\mathbb P(\tau_M^\varepsilon<T)
+
\frac{C_M(1+\mathbb E\TV(u_0))}{\Lambda},
\]
and choosing first $M$ and then $\Lambda$ yields
\eqref{eq:TV-bdd-in-prob}.
\end{proof}

\subsection{Passage to the vanishing-viscosity limit}
\label{subsec:eps-limit}

We now let $\varepsilon\downarrow0$. Since the $BV$ estimate is
localized, the viscous comparison estimate must also be localized.

For $\varepsilon,\varepsilon'>0$ and $M\ge1$, define
\[
\sigma_M^{\varepsilon,\varepsilon'}
:=
\inf
\Big\{
t\ge0:
\|u^\varepsilon(t)\|_{L^1}
+
\|u^{\varepsilon'}(t)\|_{L^1}
+
\TV(u^\varepsilon(t))
+
\TV(u^{\varepsilon'}(t))
\ge M
\Big\}
\wedge T.
\]
By
\eqref{eq:L1-untruncated-2} and
\eqref{eq:TV-bdd-in-prob},
\begin{equation}\label{eq:stop-prob-eps}
\lim_{M\uparrow\infty}
\sup_{\varepsilon,\varepsilon'>0}
\mathbb P
\big(
\sigma_M^{\varepsilon,\varepsilon'}<T
\big)
=
0.
\end{equation}

On
$[0,\sigma_M^{\varepsilon,\varepsilon'}]$,
the deterministic bounds required in
Corollary~\ref{cor:cauchy} hold. Hence
\begin{align}
&\sup_{0\le t\le T}
\mathbb E
\left\|
u^\varepsilon
(t\wedge\sigma_M^{\varepsilon,\varepsilon'})
-
u^{\varepsilon'}
(t\wedge\sigma_M^{\varepsilon,\varepsilon'})
\right\|_{L^1} \le
C_M
\left[
\mathbb E
\|u_0^\varepsilon-u_0^{\varepsilon'}\|_{L^1}
+
\sqrt{\varepsilon+\varepsilon'}
\right].
\label{eq:cauchy-stopped}
\end{align}
For the mollified initial data,
\[
\mathbb E
\|u_0^\varepsilon-u_0^{\varepsilon'}\|_{L^1}
\longrightarrow0,
\]
and under the $BV$ assumption it is in fact
$O(\varepsilon+\varepsilon')$.

Combining \eqref{eq:cauchy-stopped} with
\eqref{eq:stop-prob-eps} and the uniform second-moment estimate
\eqref{eq:L1-untruncated-2} gives
\[
u^\varepsilon
\longrightarrow
u
\quad\text{strongly in }
L^1
\big(
\Omega\times(0,T);L^1(\R)
\big)
\]
for some progressively measurable process $u$. Thus
\begin{align}
u^\varepsilon
&\longrightarrow
u
\quad\text{strongly in }
L^1
\big(
\Omega\times(0,T);L^1(\R)
\big),
\label{eq:eps-strong-L1}
\\
u^\varepsilon(t,x)
&\longrightarrow
u(t,x)
\quad
\text{for a.e. }(\omega,t,x)
\quad\text{along a subsequence}.
\label{eq:eps-ae}
\end{align}
Moreover, Fatou's lemma and
\eqref{eq:L1-untruncated-2} give
\begin{equation}\label{eq:limit-L1}
u
\in
L^2
\big(
\Omega;L^\infty(0,T;L^1(\R))
\big).
\end{equation}

\begin{proposition}[$BV$-regularity of the vanishing-viscosity limit]
\label{prop:BV-limit}
The limit process obtained above satisfies
\[
u\in L^\infty(0,T;BV(\R))
\qquad\mathbb P\text{-a.s.}
\]
\end{proposition}

\begin{proof}
Set
\[
X:=L^1\big(0,T;L^1(\R)\big)
\]
and, for $M>0$, define
\[
K_M
:=
\left\{
w\in X:
\operatorname*{ess\,sup}_{t\in(0,T)}\TV(w(t))\le M
\right\}.
\]
The set $K_M$ is closed in $X$. Indeed, if $w_n\to w$ in $X$ and
$w_n\in K_M$, then, after extraction of a subsequence,
$w_n(t)\to w(t)$ in $L^1(\R)$ for a.e. $t$. By lower semicontinuity of
total variation,
\[
\TV(w(t))
\le
\liminf_{n\to\infty}\TV(w_n(t))
\le M
\]
for a.e. $t$, and hence $w\in K_M$.

By \eqref{eq:eps-strong-L1}, the $X$-valued random variables
$u^\varepsilon$ converge to $u$ in $L^1(\Omega;X)$, hence in probability
and therefore in distribution in $X$. Since $K_M$ is closed, the
Portmanteau theorem gives
\[
\mathbb P(u\in K_M)
\ge
\limsup_{\varepsilon\downarrow0}
\mathbb P(u^\varepsilon\in K_M).
\]
On the other hand, \eqref{eq:TV-bdd-in-prob} implies that, for every
$\eta>0$, one can choose $M$ sufficiently large so that
\[
\sup_{\varepsilon>0}
\mathbb P
\left(
\sup_{0\le t\le T}\TV(u^\varepsilon(t))>M
\right)
<\eta.
\]
Consequently $\mathbb P(u\in K_M)\ge1-\eta$. Letting $M\uparrow\infty$
shows
\[
\operatorname*{ess\,sup}_{t\in(0,T)}\TV(u(t))<\infty
\qquad\mathbb P\text{-a.s.}
\]
Together with \eqref{eq:limit-L1}, this proves the assertion.
\end{proof}

\subsubsection*{Passage to the entropy inequality}

Let $\beta\in\mathcal E$, $k\in\R$, and
$0\le\phi\in C_c^\infty([0,T)\times\R)$.
Applying the It\^o--L\'evy formula to
$\beta(u^\varepsilon-k)\phi$ gives
\begin{align}
&\int_\R
\beta(u_0^\varepsilon-k)\phi(0,x)\,\d x
+
\int_{Q_T}
\beta(u^\varepsilon-k)
\partial_t\phi\,\d x\,\d t
\nonumber\\
&\quad
+
\int_{Q_T}
\mathcal U[u^\varepsilon]
\Psi_k(u^\varepsilon)
\partial_x\phi\,\d x\,\d t
\nonumber\\
&\quad
-
\int_{Q_T}
\Big[
\beta'(u^\varepsilon-k)f(u^\varepsilon)
-
\Psi_k(u^\varepsilon)
\Big]
\partial_x\mathcal U[u^\varepsilon]
\phi\,\d x\,\d t
\nonumber\\
&\quad
+
\int_0^T
\left(
\int_\R
\beta'(u^\varepsilon-k)
\Phi(u^\varepsilon)\phi\,\d x
\right)\d W(t)
\nonumber\\
&\quad
+
\frac12
\int_{Q_T}
\beta''(u^\varepsilon-k)
\Phi(u^\varepsilon)^2
\phi\,\d x\,\d t
\nonumber\\
&\quad
+
\int_0^T\!\!\int_Z
\left(
\int_\R
\Big[
\beta(
u^\varepsilon_- -k
+
G(u^\varepsilon_-,z)
)
-
\beta(u^\varepsilon_- -k)
\Big]
\phi\,\d x
\right)
\widetilde\eta(\d z,\d t)
\nonumber\\
&\quad
+
\int_0^T\!\!\int_Z\!\!\int_\R
\Big[
\beta(
u^\varepsilon-k+G(u^\varepsilon,z)
)
-
\beta(u^\varepsilon-k)
-
\beta'(u^\varepsilon-k)G(u^\varepsilon,z)
\Big]
\phi\,
\d x\,\nu(\d z)\,\d t
\nonumber\\
&\quad
+
\varepsilon
\int_{Q_T}
\beta'(u^\varepsilon-k)
\partial_{xx}u^\varepsilon
\phi\,\d x\,\d t
=
0 .
\label{eq:viscous-entropy}
\end{align}

The viscous term satisfies
\begin{align}
&\varepsilon
\int_{Q_T}
\beta'(u^\varepsilon-k)
\partial_{xx}u^\varepsilon
\phi\,\d x\,\d t
\nonumber\\
&\qquad
=
-\varepsilon
\int_{Q_T}
\beta''(u^\varepsilon-k)
|\partial_xu^\varepsilon|^2
\phi\,\d x\,\d t
-
\varepsilon
\int_{Q_T}
\beta'(u^\varepsilon-k)
\partial_xu^\varepsilon
\partial_x\phi\,\d x\,\d t .
\label{eq:viscous-entropy-decomp}
\end{align}
The first term on the right is nonpositive and therefore contributes
with the favorable entropy-dissipation sign. It need not converge to
zero.

For the second term, localize at
$\varsigma_N^\varepsilon$. By
\eqref{eq:energy-localised-eps},
\begin{align*}
&\mathbb E
\left|
\varepsilon
\int_0^{T\wedge\varsigma_N^\varepsilon}
\!\!\int_\R
\beta'(u^\varepsilon-k)
\partial_xu^\varepsilon
\partial_x\phi\,\d x\,\d t
\right| \le
C_{\beta,\phi}
\sqrt{\varepsilon}
\left(
\mathbb E
\varepsilon
\int_0^{T\wedge\varsigma_N^\varepsilon}
\|\partial_xu^\varepsilon\|_{L^2}^2\,\d t
\right)^{1/2}
\le
C_N\sqrt{\varepsilon},
\end{align*}
and this tends to zero. Since
\[
\sup_{\varepsilon}
\mathbb P(
\varsigma_N^\varepsilon<T
)
\longrightarrow0
\qquad
\text{as }N\uparrow\infty,
\]
the cross term disappears in the limit.

The strong global $L^1$ convergence also gives, by
\textbf{(K2)}--\textbf{(K3)},
\[
\|\mathcal U[u^\varepsilon]
-
\mathcal U[u]\|_{L^\infty}
\le
\|\mu\|_{L^\infty}
\|u^\varepsilon-u\|_{L^1},
\]
and
\[
\|\partial_x\mathcal U[u^\varepsilon]
-
\partial_x\mathcal U[u]\|_{L^\infty}
\le
\|\mu'\|_{L^\infty}
\|u^\varepsilon-u\|_{L^1}.
\]
Thus the two nonlocal transport contributions converge directly.

For the Brownian and jump stochastic integrals, we again localize in
the $L^1$-norm. On the localized sets, the admissible entropy
coefficients have the required uniform integrability, and
\eqref{eq:eps-strong-L1}, together with the Lipschitz assumptions on
$\Phi$ and $G$, yields convergence of the corresponding quadratic
variations. The It\^o and compensated-Poisson isometries then give
convergence of the stochastic integrals. Finally, the localization is
removed using \eqref{eq:L1-untruncated-2}.

Passing to the lower limit in
\eqref{eq:viscous-entropy} and retaining the nonnegative entropy
dissipation yields exactly the inequality
\eqref{eq:entropy-ineq}. Hence $u$ is a stochastic entropy solution.

\subsubsection*{Temporal path regularity}

It remains to verify the path requirement in
Definition~\ref{solution-defn}. The increment estimates used in
Proposition~\ref{prop:time-translation} can be repeated for
$u^\varepsilon$ after localization. Together with
\eqref{eq:L1-untruncated-2} and
\eqref{eq:energy-localised-eps}, they yield the Aldous condition in
$\mathcal S_{-1}$, uniformly in $\varepsilon$.

Moreover,
\[
\mathcal S_0=L^2(\R)
\Subset
\mathcal S_{-1},
\]
so the localized energy bounds give the corresponding compact
containment condition. Therefore the laws of
$\{u^\varepsilon\}_{\varepsilon>0}$ are tight in
\[
D([0,T];\mathcal S_{-1}).
\]
Consequently, the limit process admits an
$\mathcal S_{-1}$-valued càdlàg modification.

\subsection{Existence theorem}
\label{subsec:proof-existence}

\begin{theorem}[Existence of stochastic entropy solutions]
\label{thm:existence}
Assume
\textbf{(A1)}--\textbf{(A5)} and
\textbf{(K1)}--\textbf{(K4)}.
Then the stochastic nonlocal conservation law
\eqref{eq:u}--\eqref{eq:id} admits a vanishing-viscosity stochastic
entropy solution $u$ in the sense of
Definition~\ref{solution-defn}. Moreover,
\begin{equation}\label{eq:existence-L1}
u
\in
\mathcal N_\Omega^2(0,T;L^1(\R))
\cap
L^2
\big(
\Omega;L^\infty(0,T;L^1(\R))
\big), \quad \mbox{and} \quad u
\in
D([0,T];\mathcal S_{-1})
\qquad\mathbb P\text{-a.s.}
\end{equation}
Moreover,
\[
u\in L^\infty(0,T;BV(\R))
\qquad\mathbb P\text{-a.s.}
\]
The viscous approximations satisfy
\[
u^\varepsilon
\longrightarrow
u
\quad\text{strongly in }
L^1
\big(
\Omega\times(0,T);L^1(\R)
\big).
\]
Furthermore, their spatial variations are uniformly bounded in
probability in the sense of \eqref{eq:TV-bdd-in-prob}.
\end{theorem}

\begin{proof}
For every fixed $\varepsilon>0$, the stopped Galerkin construction of
Sections~\ref{subsec:wellposed}--\ref{subsec:galerkin-compactness}
produces a maximal viscous solution. Proposition~\ref{prop:L1}
provides an $L^1$ estimate independent of the $L^2$-exit level, and
the localized energy estimate then rules out finite-time explosion.
Thus Proposition~\ref{prop:remove-truncation} yields a unique global
viscous solution $u^\varepsilon$.

Under \textbf{(K4)}, Proposition~\ref{prop:BV} gives the localized
$BV$-estimate and uniform control of the corresponding exit probabilities.
Combining this estimate with the stopped
$L^1$-comparison of Corollary~\ref{cor:cauchy} shows that
$\{u^\varepsilon\}$ is Cauchy in
\[
L^1
\big(
\Omega\times(0,T);L^1(\R)
\big).
\]
Hence, $u^\varepsilon\to u$ strongly in this space.

The viscous entropy identity passes to the limit as described above;
the entropy-dissipation contribution has the favorable sign and the
remaining viscous cross term converges to zero. Therefore $u$
satisfies \eqref{eq:entropy-ineq}. The second-moment $L^1$ estimate and
Fatou's lemma yield \eqref{eq:existence-L1}. Proposition~\ref{prop:BV-limit}
shows in addition that $u\in L^\infty(0,T;BV(\R))$ almost surely, while
the localized Aldous estimates provide an $\mathcal S_{-1}$-valued càdlàg
modification.

Thus $u$ is a stochastic entropy solution of
\eqref{eq:u}--\eqref{eq:id}.
\end{proof}

\medskip

\section{Concluding Remarks and Future Directions}
\label{sec:conclusion}
We studied a one-dimensional stochastic nonlocal conservation law driven
by multiplicative Brownian and L\'evy jump noise. We introduced an entropy
formulation adapted to the solution-dependent nonlocal velocity, based on
the shifted entropy flux $\Psi_k$ and an additional
$\partial_x\mathcal U[u]$-term, which is essential for the $L^1$-stability
analysis. Existence of spatially $BV$ stochastic entropy solutions was obtained via
vanishing viscosity, and uniqueness in the class of vanishing-viscosity
solutions via a spatial doubling argument at equal times, which avoids
the adaptedness issues of space--time doubling. Under \textbf{(A5$'$)},
the same estimate yields the localized Cauchy rate
$O(\sqrt{\varepsilon+\varepsilon'})$ for the viscous approximations, and
hence a localized convergence rate of order $\sqrt{\varepsilon}$ in the
vanishing-viscosity limit.

Several questions remain open:

\begin{enumerate}[label=\textup{(\arabic*)}, leftmargin=2.4em]

\item \textbf{Uniqueness in the full entropy class.} Extending the
comparison principle to arbitrary entropy solutions of
Definition~\ref{solution-defn} would require a direct stochastic doubling
argument with an adapted random entropy level, as in
\cite{BVW2012,BVW2014}, together with a treatment of the nonlocal term.

\item \textbf{Weaker assumptions.} It would be of interest to obtain a
fully localized quantitative stability estimate under moment assumptions
only, removing \textbf{(A5$'$)}, and to relax the kernel regularity in
\textbf{(K3)}--\textbf{(K4)}, e.g.\ to $BV$ kernels.

\item \textbf{Nonlocal-to-local limit and numerics.} A natural problem is
the behavior of solutions as the kernel concentrates to a Dirac mass, in
analogy with the deterministic theory
\cite{CCS2019,keimer2019approximation,CGES2021,bressan2021entropy,
Col23a,CCDKP2024,CCMS2024}, as well as the convergence analysis of
finite-volume schemes for the stochastic nonlocal model.

\end{enumerate}

\section*{Statements and Declarations}

\subsection*{Funding}
The second author gratefully acknowledges the financial support provided by the Anusandhan
National Research Foundation (ANRF) under the Early Career Research Grant (ECRG), Grant
No.~ANRF\slash ECRG\slash 2025\slash 002863\slash PMS. The third author gratefully acknowledges the support of the DST-INSPIRE Faculty Research Grant
(Registration No.~IFA21-MA158), Department of Science and Technology, Government of India.

\subsection*{Competing interests}
The authors declare that they have no competing interests, financial or non-financial, that are relevant to the content of this article.


\medskip

\bibliographystyle{siam} \bibliography{jcp,stoch}

\end{document}